\documentclass[12pt]{article} 

\usepackage{amsmath}
\usepackage{amsthm}
\usepackage{amsfonts}
\usepackage{mathrsfs}
\usepackage{stmaryrd}
\usepackage{setspace}
\usepackage{fullpage}
\usepackage{amssymb}
\usepackage{breqn}
\usepackage{enumitem}
\usepackage{bbold} 
\usepackage{authblk}
\usepackage{comment}
\usepackage{hyperref}
\usepackage{pgf,tikz}
\usepackage{graphicx}
\usepackage{xcolor}
\usepackage{caption}
\usetikzlibrary{decorations.pathreplacing,arrows}
\usetikzlibrary{positioning,chains,fit,shapes,calc}
\tikzstyle{vertex}=[circle,draw=black,fill=black,inner sep=0,minimum size=3pt,text=white,font=\footnotesize]
\usetikzlibrary{graphs}
\usetikzlibrary{shapes.symbols}

\newtheorem*{rep@theorem}{\rep@title}
\newcommand{\newreptheorem}[2]{%
\newenvironment{rep#1}[1]{%
 \def\rep@title{#2 \ref{##1}}%
 \begin{rep@theorem}}%
 {\end{rep@theorem}}}
\makeatother

\theoremstyle{plain}

\newtheorem*{proposition*}{Proposition}

\newtheorem{thm}{Theorem}[section]
\newtheorem{theorem}{Theorem}[section]

\newtheorem{lemma}[thm]{Lemma}
\newtheorem{conj}[thm]{Conjecture}

\newtheorem{prop}[thm]{Proposition}

\newtheorem{defn}[thm]{Definition}

\newtheorem{claim}[thm]{Claim}

\newtheorem*{theorem*}{Theorem}
\newtheorem*{prop*}{Proposition}

\newcommand\ex{\ensuremath{\mathrm{ex}}}

\newcommand\cF{{\mathcal F}}
\newcommand\cG{{\mathcal G}}

\newcommand\cJ{{\mathcal J}}

\newcommand\cN{{\mathcal N}}
\newcommand\cP{{\mathcal P}}

\newcommand\cT{{\mathcal T}}

\def\lf{\left\lfloor}   
\def\rf{\right\rfloor}

\newcommand{\abs}[1]{\left\lvert{#1}\right\rvert}

\newcommand{\ignore}[1]{}

\title{Generalized planar Tur\'an numbers related to short cycles}
\author{Ervin Gy\H{o}ri\footnote{Alfr\'ed R\'enyi Institute of Mathematics, HUN-REN, Budapest, Hungary E-mail: \texttt{gyori@renyi.hu}}, Hilal Hama Karim\footnote{Department of Computer Science and Information Theory, Faculty of Electrical Engineering and Informatics, Budapest University of Technology and Economics, Műegyetem rkp. 3., H-1111 Budapest, Hungary. E-mail: \texttt{hilal.hamakarim@edu.bme.hu}}}

\date{}

\begin{document}

\maketitle

 \begin{abstract}
    Given two graphs $H$ and $F$, the generalized planar Tur\'an number $\ex_\cP(n,H,F)$ is the maximum number of copies of $H$ that an $n$-vertex $F$-free planar graph can have. We investigate this function when $H$ and $F$ are short cycles. Namely, for large $n$, we find the exact value of $\ex_\cP(n, C_l,C_3)$, where $C_l$ is a cycle of length $l$, for $4\leq l\leq 6$, and determine the extremal graphs in each case. Also, considering the converse of these problems, we determine sharp upper bounds for $\ex_\cP(n,C_3,C_l)$, for $4\leq l\leq 6$.

\end{abstract}

\section{Introduction}
Tur\'an type problems are fundamental in extremal graph theory. Given a family $\cF$ of graphs, A graph $G$ is said to be $\cF$-free if it does not contain any member of $\cF$ as a subgraph. The ordinary Tur\'an number, $\ex(n,\cF)$ (simply, \ex(n,F), if $\cF=\{F\}$), is the maximum number of edges an $\cF$-free simple graph on $n$ vertices can have. This has been of great interest ever since its start from the renowned theorem of Tur\'an \cite{tur41} in 1941, who determined $\ex(n,K_r)$, where $K_r$ is a complete graph on $r$ vertices (see \cite{fursim} for a survey). A natural generalization is to consider maximizing the number of copies of another subgraph $H$ instead of edges while forbidding $\cF$, which is denoted by $\ex(n,H,\cF)$. Before initiating its general and systematic study by Alon and Shikhelman \cite{alon} in 2016, there had been sporadic results regarding this function. For instance, Erd\H{o}s \cite{Erd62}, determined $\ex(n,K_r,K_t)$ for all $r<t$, and Gy\H{o}ri, Pach and Simonovits \cite{GyPS} studied $\ex(n,H,K_r)$, especially when $r=3$ and $H$ is bipartite. Probably the most interesting problems among all of these ones are those concerning pentagons $C_5$ (cycles of length 5) and triangles $C_3$. In 1984, Erd\H{o}s \cite{Erd84}  conjectured that the maximum number of $C_5$'s in a triangle free graph on $n$ vertices is at most $(n/5)^5$, which took several attempts until it was proved in 2013, see \cite{Gyo, grz, hhknr}. Bollob\'as and Gy\H{o}ri \cite{BoGy} studied the converse of this problem, $\ex(n,C_3,C_5)$.

The study of the planar version of Tur\'an numbers was initiated by Dowden \cite{Dowden} in 2015. Namely, he studied the function $\ex_\cP(n,\cF)$, the maximum number of edges in $\cF$-free $n$-vertex planar graphs. He obtained upper bounds for $\ex_\cP(n,C_4)$ and $\ex_\cP(n,C_5)$ with constructions attaining the bounds for infinitely many $n$. Ghosh, Gy\H{o}ri, Martin, Paulos and Xiao \cite{GhGyMPX} determined $\ex_\cP(n,C_6)$ introducing some contribution method. Applying this method, Gy\H{o}ri, Wang and Zheng \cite{GyWZh} obtained various other results and gave a new shorter proof for $\ex_\cP(n,C_5)$ as well as a clearer extremal construction.

More recently, Gy\H{o}ri, Paulos, Salia, Tompkins and Zamora \cite{Gypntz} introduced the generalized version of the planar Tur\'an numbers, $\ex_\cP(n,H,\cF)$, the maximum number of copies of a subgraph $H$ in $\cF$-free $n$-vertex planar graphs. In particular, they showed that for any $k\geq 5$, $\ex_\cP(n,C_k, C_4)=\Theta (n^{\lfloor k/3\rfloor})$, and in case $k=5$, they proved $\ex_\cP(n,C_5,C_4)=n-4$, for all $n \geq 5$ (except $n=6$).  Maximizing the number of copies of a given subgraph $H$ in planar graphs can be viewed as a special case of generalized planar Tur\'an problems by taking $\cF= \emptyset$. Hakimi and Schmeichel \cite{Hakimi} determined the maximum number of triangles and 4-cycles in planar graphs. Gy\H{o}ri, Paulos, Salia, Tompkins and Zamora \cite{Gypntz1} determined the maximum number of $5$-cycles in planar graphs, proving that $\ex_\cP(n,C_5, \emptyset)=2n^2-10n+12$, for every $n=6$ or $n\geq 8$. While exact results for longer cycles are not known, Cox and Martin \cite{coxmartin} developed a general technique to count subgraphs in planar graphs and they conjectured that the maximum number of an even cycle $C_{2k}$ is asymptotically $(n/k)^k$. This conjecture was proved by Lv, Győri, He, Salia, Tompkins and Zhu \cite{LGHSTZ}.

\begin{theorem}\label{asymptC2k}\cite{LGHSTZ} For every $k\geq 3$, $\ex_\cP(n,C_{2k},\emptyset)=\left(\frac{n}{k}\right)^k+o(n^k)$.
    \end{theorem}

Another direction of research about the induced subgraphs in planar graphs. Ghosh, Gy\H{o}ri, Janzer, Paulos, Salia, Zamora \cite{GhGyJPSZ}, and independently Savery \cite{savery2}, determined the maximum number of induced 5-cycles in planar graphs. Savery \cite{savery} extended this to induced 6-cycles.  

Here, we consider triangle-free planar graphs and count cycles of lengths 4, 5 or 6. Similarly, we determine the maximum possible number of triangles in planar graphs when such a cycle is forbidden.

\begin{theorem}\label{C4C3}
    For every $n \geq 4$, $\displaystyle \ex_\cP(n,C_4,C_3)= \binom{n-2}{2}$, and the unique extremal graph is $K_{2,n-2}$. 
\end{theorem}

\begin{theorem}\label{C4C5}
    For every $n \geq 4$, $\displaystyle \ex_p(n,C_4,C_5)= \binom{n-2}{2}$. 
\end{theorem}

\begin{theorem}\label{C3C4}
    For every $n \geq 4$, we have $\ex_\cP(n,C_3,C_4) \leq \frac{5}{7}(n-2)$, and this bound is sharp for infinitely many values of $n$.
\end{theorem}

The following theorem is a bit different, as we determine the maximum number of triangles in $K_4$-free planar graphs.

\begin{theorem}\label{C3K4}
    For every $n \geq 3$, we have $\ex_\cP(n,C_3,K_4) \leq \frac{7}{3}n -6$, and this bound is sharp for all $n$ divisible by $3$.
\end{theorem}

\begin{theorem}\label{C5C3weak}
    For every $n \geq 5$, $\displaystyle \ex_p(n,C_5,C_3) = \lf (n-3)/2 \rf \cdot \lceil (n-3)/2 \rceil$.
\end{theorem}
An extremal graph for $\ex_p(n,C_5,C_3)$ can easily be seen to be a $C_5$ on which two non-adjacent vertices are blown-up in a balanced way. Somewhat surprisingly, there are many other extremal graphs.

For each $n\geq 5$, define a class $\cJ_n$ of $n$-vertex planar graphs $J_n$ as follows. Take a regular pentagon, replace one of its vertices, $x$, by an independent set of vertices $C$, and replace the edge $yz$ opposite to $x$ by a tree with color classes $A$ and $B$, such that $|C|$ and $|A\cup B|-1$ are as equal as possible (see Figure \ref{EXTC6C5}(a)). We prove the following stronger theorem.
\begin{theorem}\label{C5C3}
For every $n \geq 5$, $\displaystyle \ex_\cP(n,C_5,C_3) = \lf (n-3)/2 \rf \cdot \lceil (n-3)/2 \rceil$, and  $\cJ_n$ is the set of all extremal graphs.    
\end{theorem}


In the other direction, maximizing the number of triangles  pentagon-free graphs, we prove the following theorem.
\begin{theorem}\label{C3C5}
    For every $n\geq 11$, $\ex_\cP(n,C_3,C_5)\leq \lf \frac{8n-22}{5} \rf$, and this bound is sharp for infinitely many values of $n$.
\end{theorem}

Then, we turn to triangles and 6-cycles, and first maximize the number of triangles in $C_6$-free graphs.

\begin{theorem}\label{C3C6}
    For every $n\geq 18$, $\ex_\cP(n,C_3,C_6) \leq \frac{35n-98}{18}$, and this bound is sharp for infinitely many values of $n$.
\end{theorem}
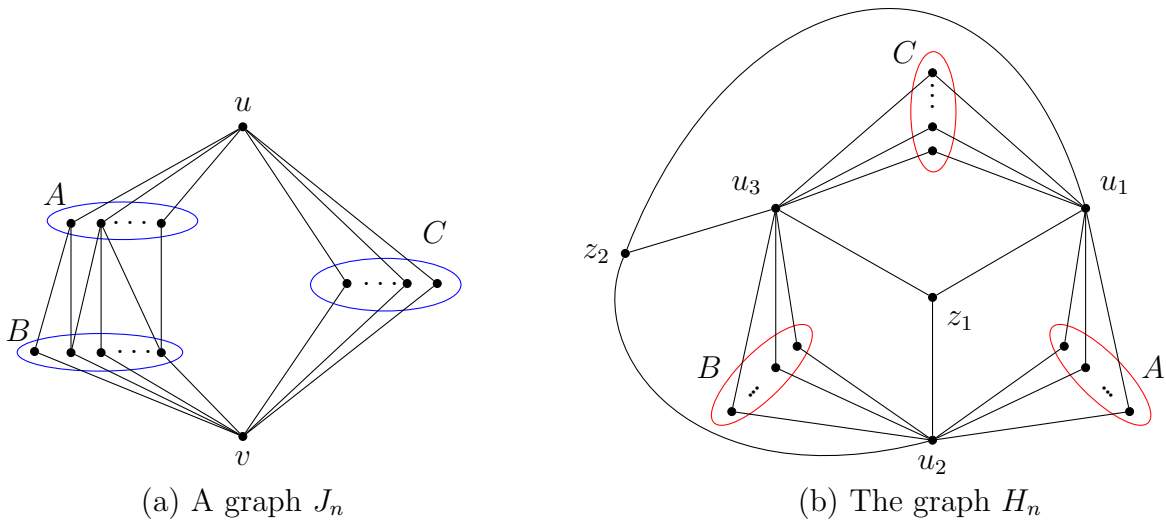
\begin{figure}[ht]
\begin{center}
    \begin{tikzpicture}
    \node[vertex, label=above: $u$] (u) {};
\node [vertex] (c1) [below right= 2 and 1.3cm of u] {};
\node [] (c2) [below right= 1.87 and 1.4cm of u] {$\ldots$};
\node [vertex] (c3) [below right= 2 and 2.1cm of u] {};
\node [vertex] (c4) [below right= 2 and 2.5cm of u] {};

\node [vertex, label=below: $v$] (v) [below =4cm of u] {};
\node [vertex] (a1) [below left=1.2cm and 1cm of u] {};
\node [] (a2) [below left=1.07cm and 1.07cm of u] {$\ldots$};
\node [vertex] (a3) [below left=1.2cm and 1.8cm of u] {};
\node [vertex] (a4) [below left=1.2cm and 2.2cm of u] {};

\node [vertex] (b1) [below=1.6cm of a1] {};
\node [] (b2) [below right= 1.5 and 0.01cm of a3] {$\ldots$};
\node [vertex] (b3) [below=1.6cm of a3] {};
\node [vertex] (b4) [below=1.6cm of a4] {};
\node [vertex] (b5) [below left=1.62cm and 0.4cm of a4] {};

\node[vertex, label=above left: $u_3$] (u3) [below right=1cm and 7cm of u] {};
\node[vertex] (v3) [above right=1cm and 2cm of u3] {};
\node[vertex] (v31) [above=0.6cm of v3] {};
\node[] (v32) [above=0.05cm of v3] {$\vdots$};
\node[vertex] (v33) [below=0.2cm of v3] {};

\node[vertex, label=above right: $u_1$] (u1) [right=4cm of u3] {};
\node[vertex, label=below: $u_2$] (u2) [below right=3cm and 2cm of u3] {};

\node[vertex, label=below right: $z_1$] (z1) [below right=1.1cm and 2cm of u3] {};
\node[vertex, label=left: $z_2$] (z2) [above left=0.5cm and 4cm of z1] {};

\node[vertex] (v2) [below=2cm of u3] {};
\node[vertex] (v21) [above right=0.2cm and 0.2cm of v2] {};
\node[] (v22) [below left=0.001cm and 0.001cm of v2] {$\cdot$};
\node[] (v222) [below left=0.05cm and 0.05cm of v2] {$\cdot$};
\node[] (v223) [below left=0.09cm and 0.09cm of v2] {$\cdot$};
\node[vertex] (v23) [below left=0.5cm and 0.5 of v2] {};

\node[vertex] (v1) [below=2cm of u1] {};
\node[vertex] (v11) [above left=0.2cm and 0.2cm of v1] {};
\node[] (v12) [below right=0.001cm and 0.001cm of v1] {$\cdot$};
\node[] (v122) [below right=0.05cm and 0.05cm of v1] {$\cdot$};
\node[] (v123) [below right=0.09cm and 0.09cm of v1] {$\cdot$};
\node[vertex] (v13) [below right=0.5cm and 0.5 of v1] {};

 \node[] at (2.55,-1.4) {$C$};
 \node[] at (-2.5,-0.9) {$A$};
 \node[] at (-2.99,-2.7) {$B$};

\draw[blue] (1.9,-2.1) ellipse (1cm and 0.35cm);
\draw[blue] (-1.6,-1.25) ellipse (1cm and 0.25cm);
\draw[blue] (-1.9,-3) ellipse (1.1cm and 0.25cm);

\draw[color=red] (9.18,0.2) ellipse (0.3cm and 0.8cm);
\draw[rotate around={45:(6.9,-3.3)},red] (6.9,-3.3) ellipse (0.9cm and 0.3cm);
\draw[rotate around={135:(11.4,-3.3)},red] (11.4,-3.3) ellipse (0.9cm and 0.3cm);

\node[] at (12.1,-3.2) {$A$};
\node[] at (6.2,-3.2) {$B$};
\node[] at (8.8,1) {$C$};

\path (u) edge (c1);
\path (u) edge (c3);
\path (u) edge (c4);
\path (u) edge (a1);
\path (u) edge (a3);
\path (u) edge (a4);
\path (v) edge (c1);
\path (v) edge (c3);
\path (v) edge (c4);
\path (v) edge (b1);
\path (v) edge (b3);
\path (v) edge (b4);
\path (v) edge (b5);
\path (b1) edge (a1);
\path (b1) edge (a3);
\path (a3) edge (b3);
\path (a3) edge (b4);
\path (b4) edge (a4);
\path (a4) edge (b5);

\path (u1) edge (z1);
\path (u2) edge (z1);
\path (u3) edge (z1);
\path (u1) edge (v1);
\path (u1) edge (v11);
\path (u1) edge (v13);
\path (u1) edge (v3);
\path (u1) edge (v31);
\path (u1) edge (v33);
\path (u3) edge (v3);
\path (u3) edge (v31);
\path (u3) edge (v33);
\path (u2) edge (v1);
\path (u2) edge (v11);
\path (u2) edge (v13);
\path (u3) edge (v2);
\path (u3) edge (v21);
\path (u3) edge (v23);
\path (u2) edge (v2);
\path (u2) edge (v21);
\path (u2) edge (v23);
\path (u3) edge (z2);
\draw (u2) .. controls (6,-5) and (4.5,-3) .. (z2);
\draw (u1) .. controls (10,3) and (6.5,2) .. (z2);

\node[] at (0,-5) {(a) A graph $J_n$};
\node[] at (9,-5) {(b) The graph $H_n$};
\end{tikzpicture} 
\caption{The extremal graphs for $\ex_\cP(n,C_5,C_3)$ and $\ex_\cP(n,C_6,C_3)$}
\label{EXTC6C5}
\end{center} 
\end{figure}

Finally, we determine the maximum number of 6-cycles while forbidding triangles, together with the unique extremal graph achieving the maximum value. Along the proof, we will be considering the number of paths of length four (i.e. of five vertices) between two vertices. This is an interesting problem on its own, and we prove a theorem as follows.

\begin{theorem}\label{P5inC3-free}
    Let $G$ be a triangle-free planar graph on $n\geq 5$ vertices. For any two vertices $u, v \in V(G)$, there are at most $(\frac{n-1}{2})^2-2$ paths of length four connecting them.
\end{theorem}

Before the stating the next theorem, we present the following construction. For every $n\geq 6$, define a graph $H_n$ as follows. The vertex set of $H_n$ consists of $\{u_1,u_2,u_3\} \cup \{z_1,z_2\} \cup A \cup B \cup C$, such that each of these sets is an independent set of vertices, they are pairwise disjoint, each of the $u_i$'s is adjacent to each of the $z_i$'s, every vertex in $A$ is adjacent to both of $u_1$ and $u_2$, every vertex in $B$ is adjacent to both of $u_2$ and $u_3$, and every vertex in $C$ is adjacent to both $u_1$ and $u_3$. Moreover, the sizes of $A$, $B$ and $C$ are as equal as possible (see Figure  \ref{EXTC6C5}(b)). It is easy to see that $H_n$ contains $h(n)$ 6-cycles (nevertheless, we will give an explanation for this in Section \ref{sec4}), where $h(n)$ is defined as follows.
\begin{align*}\displaystyle
    h(n)=\left\{
\begin{array}{ll}
    \frac{n^3}{27}+\frac{n^2}{9}-2n+2, &~ \ \ \text{if} \ n \equiv 0 \ (mod \ 3) \\
     \frac{n^3}{27}+\frac{n^2}{9}-2n+\frac{50}{27}, &~ \ \ \text{if} \ n \equiv 1 \ (mod \ 3)\\
      \frac{n^3}{27}+\frac{n^2}{9}-\frac{17n}{9}+\frac{55}{27}, &~ \ \ \text{if} \ n \equiv 2 \ (mod \ 3)
\end{array} 
\right.
\end{align*}

\begin{theorem}\label{C6C3}
    Let $n$ be sufficiently large. Then, $\ex_\cP(n,C_6,C_3)=h(n)$, and the unique extremal graph is $H_n$. 
\end{theorem}

We think that theorems \ref{C5C3} and \ref{C6C3} are the highlights of this paper. The proof of theorem \ref{C6C3} makes it possible to prove the following interesting version of Theorem \ref{P5inC3-free}.

\begin{theorem}\label{P53vtxC3free}
    Let $G$ be a triangle-free planar graph on $n$ vertices. For any three distinct vertices $u_1, u_2, u_3 \in V(G)$, the number of paths of length four joining all the three pairs of them is at most $3\left(\frac{n+1}{3}\right)^2-6$. Moreover, for each $n\equiv 2$ (mod 3), the graph $H_n$ attains this bound.
\end{theorem}

Throughout this paper we use the following notations. $V(G)$ and $E(G)$ denote the vertex and edge sets of a graph $G$, respectively, and $e(G):=|E(G)|$ (sometimes just $e$). For a subset $X\subseteq V(G)$, we denote by $G[X]$ the subgraph of $G$ induced  on $X$, and $G\setminus X$ (or simply $G-v$, if $X=\{v\}$) denotes the induced subgraph $G[V(G)\setminus X]$. $N(x)$ denotes the set of the neighbors of a vertex $x$, and $N(x_1, \ldots, x_k)$ denotes the common neighbors of all the vertices $x_1, \ldots ,x_k$. For subsets $A, B \subseteq V(G)$, the edges of $G$ between vertices of $A$ and $B$ is denoted by $E(A,B)$. Given graphs $H$ and $G$, we denote the number of copies of $H$ in $G$ (i.e. subgraphs of $G$ isomorphic to $H$) by $\cN(H,G)$, when the graph $G$ is clear we also use $\# H$ to denote the number of copies of $H$. For a positive integer $t$, we use $[t]$ to denote the set $\{1,2,\ldots,t\}$. We say that a cycle $C$ in a plane graph $G$ is not empty if there are vertices of $G$ in the region bounded by $C$, and $C$ is said to be a separating cycle if there are vertices in both the interior and the exterior regions of $C$. We let $P_l$ denote a path on $l$ vertices (of length $l-1$), and $P_l(x, y)$ denote a path $P_l$ connecting the two vertices $x$ and $y$. A cycle of length $l$ is denoted by $C_l$. Given a vertex $v$, we denote by $C_l(v)$ a cycle of length $l$ that contains $v$. 

\section{4-cycles}\label{sec2}

In this section we prove Theorems \ref{C4C3} through \ref{C3K4}.

\begin{proof}[Proof of Theorem \ref{C4C3}]
    Obviously, the complete bipartite graph $K_{2,n-2}$ is a planar triangle free graph that contains $\binom{n-2}{2}$ copies of $C_4$. Let $G$ be a triangle-free planar graph on $n\geq 4$ vertices. We prove that $\cN(C_4,G) \leq \binom{n-2}{2}$, and equality holds if and only if $G$ is $K_{2,n-2}$.
    
    The proof is by induction on $n$. For $n \leq 5$ the statement holds obviously. Assume $n > 5$. Let $G_0$ be a maximum $K_{2,k}$ subgraph of $G$. 
    Let $\{u,v\}$ be the class of $G_0$ of size two, and $\{x_1,x_2, \ldots, x_k\}$ be its class of size $k$. Since $G$ is triangle-free, these sets are independent in $G$, too. 
    Note that the subgraph $G_0$ has $k$ faces each of which is a $C_4$. 
    If a face of $G_0$ is not a face in $G$, i.e. it contains some vertices of $G$ in its interior, then there cannot be any crossing $C_4$,i.e. a $4$-cycle that uses vertices from both the interior and the exterior of the face. Indeed such a crossing $C_4$ must use exactly one vertex within the face, and since $G$ is planar and triangle-free, it must be of the form $yuzv$, for some vertex $y$ within the face and $z$ not in the face. Then, $y$ is a common neighbor of $u$ and $v$ different from $x_i$'s, contradicting the maximality of $G_0$. Note that if all the 4-cycles of $G$ are faces of $G$, then due to Euler's formula, $\cN(C_4,G)\leq n-2 < \binom{n-2}{2}$ (since $n>5$), so we may assume that $G$ contains 4-cycles that are not faces. 

    \begin{claim}\label{clmsepC_4}
        If $V(G) \setminus V(G_0) \neq \emptyset$, then $G$ has a separating $C_4$ for which there are no crossing $4$-cycles.
    \end{claim}
    \begin{proof}[Proof of the claim.] If $k=2$, then $G_0$ is a $C_4$ and has two faces. We may choose it to be a separating $C_4$, i.e. both faces are non-empty, since otherwise all the 4-cycles of $G$ are faces, a contradiction. 
    
    If $k>2$, then it has a non-empty face, since $V(G) \setminus V(G_0) \neq \emptyset$. Hence, obviously that face is a separating $C_4$. 
    In both cases the separating $C_4$ is a face of $G_0$, and hence as explained above, for this separating $C_4$ there are no crossing $4$-cycles in $G$.
        \end{proof}
    Now, let us complete the proof of the theorem. If $V(G) \setminus V(G_0) = \emptyset$, then $k=n-2$, $G$ is the graph $K_{2,n-2}$ and $\cN(C_4,G)=\binom{n-2}{2}$. Otherwise, $G$ has a separating $C_4$ by Claim \ref{clmsepC_4}, call it $C$, for which there are no crossing $4$-cycles in $G$. Let $L$ be the set of vertices of $G$ in the interior of $C$, with $|L|=l$. Let $G_1:=G[L \cup V(C)]$ and $G_2:=G[V(G)\setminus L]$. Thus, 
    \begin{align*}\label{ineq1}
        \cN(C_4,G)&= \cN(C_4, G_1)+\cN(C_4,G_2) -1 \\
        &\leq \binom{(l+4)-2}{2} + \binom{(n-l)-2}{2} -1 \\
        &= \binom{l+2}{2} + \binom{(n-2)-l}{2} -1  \\
        &=\binom{n-2}{2}+l^2+4l-nl \\
        &< \binom{n-2}{2}.
   \end{align*}
   where the term  -1 is to avoid double counting $C$, the first inequality is due to the induction hypothesis, and the last one is justified by the fact that $n\geq l+5$ and $l \geq 1$ (since $C$ is a separating cycle), which gives $l^2+4l-nl \leq l^2+4l-(l+5)l = -l < 0$.
\end{proof}

\begin{proof}[Proof of Theorem \ref{C4C5}] The graph $K_{2,n-2}$ is a planar $C_5$-free graph that contains $\binom{n-2}{2}$ 4-cycles, showing that $\ex_\cP(n,C_4,C_5) \geq \binom{n-2}{2}$. Let $G$ be a $C_5$-free planar graph on $n\geq 4$ vertices. Observe that a triangle and a $C_4$ cannot share exactly one edge, since otherwise they would worm a $C_5$. Obviously, a triangle does not share all its three edges with a $C_4$. Then, for any triangle in $G$, if each of its edges is in some $C_4$, then it must share exactly one edge with some $C_4$, forming a $C_5$. Consequently, every triangle of $G$ has at least one edge that is not contained in any $C_4$. Deleting such edges from each triangle, we can make $G$ triangle-free without reducing the number of 4-cycles. Thus, $\ex_\cP(n,C_4,C_5) \leq \ex_\cP(n,C_4,C_3)=\binom{n-2}{2},$
    where the last equality comes from Theorem \ref{C4C3}.
\end{proof}

\begin{proof}[Proof of Theorem \ref{C3C4}]
    If $G$ is a planar $C_4$-free graph, then no two triangles in $G$ can share an edge, as otherwise a $C_4$ would be formed. Then, every edge of $G$ is in at most one triangle, which gives $$3\cN(C_3,G)\leq e(G) \leq \ex_\cP(n, C_4).$$ Dowden \cite{Dowden} showed that $\ex_\cP(n,C_4)\leq \frac{15}{7}(n-2)$, which implies $\cN(C_3,G)\leq \frac{5}{7}(n-2)$, and hence,
    \[\ex_\cP(n,C_3,C_4) \leq \frac{5}{7}(n-2).\]
    They also provided constructions for every $n\equiv 30 \ (mod \ 70)$ that attain $\ex_\cP(n,C_4)$, implying the sharpness of their result. In fact, every edge in their extremal graphs is in a $C_3$ (see Theorem(2) in \cite{Dowden}), and hence they provide graphs that achieve the stated upper bound in our case as well, implying the sharpness of the bound.
\end{proof}

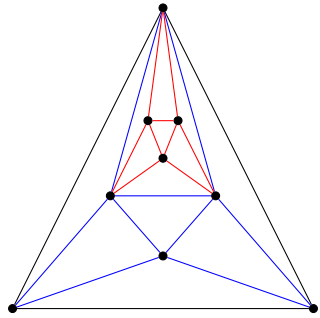
\begin{figure}[ht]
\begin{center}
    \begin{tikzpicture}
    \node[vertex] (a) at (0,4) {};
    \node[vertex] (b) at (-2,0) {};
    \node[vertex] (c) at (2,0) {};

    \node[vertex] (a1) at (-0.7,1.5) {};
    \node[vertex] (b1) at (0.7,1.5) {};
    \node[vertex] (c1) at (0,0.7) {};

    \node[vertex] (a2) at (-0.2,2.5) {};
    \node[vertex] (b2) at (0.2,2.5) {};
    \node[vertex] (c2) at (0,2) {};

    \draw (a) -- (b) -- (c) -- (a);
    \draw[blue] (a1) -- (b1) -- (c1) -- (a1);
    \draw[blue] (a1) -- (a);
    \draw[blue] (a1) -- (b);
    \draw[blue] (b1) -- (a);
    \draw[blue] (b1) -- (c);
    \draw[blue] (c1) -- (c);
    \draw[blue] (c1) -- (b);

    \draw[red] (a2) -- (b2) -- (c2) -- (a2);
    \draw[red] (a) -- (a2);
    \draw[red] (a) -- (b2);
    \draw[red] (a1) -- (a2);
    \draw[red] (a1) -- (c2);
    \draw[red] (b1) -- (b2);
    \draw[red] (b1) -- (c2);
    \end{tikzpicture}
    \end{center}
    \caption{A $K_4$-free graph attaining maximum number of triangles.} 
    \label{EXTC3K4}
    \end{figure}
    
\begin{proof}[Proof of Theorem \ref{C3K4}]
    The proof is by induction on $n$. Let $G$ be a $K_4$-free planar graph on $n$ vertices. We can easily check that the result holds for $n=3, 4$ and $5$, so assume $n \geq 6$. Further, we may assume that $G$ is connected as otherwise we may apply induction on the components. If every $C_3$ in $G$ is a face, then using Euler's formula, we have $\cN(C_3,G)=f=e+2-n \leq (3n-6)+2-n\leq 7n/3-6$. Let $C$ be a triangle that is not a face, then it contains vertices both in its interior and its exterior. Let $C$ contain $m$ vertices inside and $n-m$ vertices outside. Let $G_1$ be the induced subgraph of $G$ on the vertices of $C$ and those in its interior, and $G_2$ be the induced subgraph on the vertices on $C$ and those in its exterior. Note that $|V(G_1)|=m+3, |V(G_2)|=n-m$ and $C$ is the only triangle that is in both of $G_1$ and $G_2$. 
    Applying the induction hypothesis, we obtain 
    \begin{align*}
        \cN(C_3,G)&=\cN(C_3,G_1)+\cN(C_3,G_2)-1 \\
        &\leq \left( \frac{7}{3}(m+3)-6 \right)+ \left(\frac{7}{3}(n-m)-6\right)-1\\
        &=\frac{7}{3}n-6.
    \end{align*}

    The sharpness of the bound is shown by the following construction. Starting from a triangle, insert a triangle inside and in the region between the two triangles join all the vertices via a 6-cycle, keeping the graph planar. Now, we obtain $G_{k+1}$ from $G_k$ by inserting a triangle in a face of $G_k$ (all of which are triangles) in the same way, as shown in  Figure \ref{EXTC3K4}. 
    \end{proof}

\section{Triangles and 5-cycles}\label{sec3}

In this section we prove Theorem \ref{C5C3} and Theorem \ref{C3C5}. For the first one we follow a similar idea as in the proof of Theorem \ref{C4C3}. Basically, given an $n$-vertex graph $G$ we take a maximal subgraph $G_0$ that is a slight modification of a member of $\cJ_n$, and then prove that $G$ has the maximum number of $C_5$'s if $G_0$ contains all vertices of $G$.

\begin{proof}[Proof of Theorem \ref{C5C3}] 
The proof is by induction on $n$. Let $G$ be a triangle free plane graph on $n$ vertices. The result is obviously true for $n=5$, so assume $n \geq 6$. Let $u$ and $v$ be two non-consecutive vertices on a $C_5$ in $G$, and then take $G_0$ to be the subgraph of $G$ that consists of all the $5$-cycles of $G$ that contain $u$ and $v$. Note that since $G$ is triangle-free, then $u$ and $v$ are not adjacent, and $G_0$ consists of all the paths of length 2 and 3 connecting $u$ and $v$ in $G$. Observe also that a common neighbor of $u$ and $v$ cannot be on a path of length $3$ joining them, as otherwise a triangle would be created. Therefore, $G_0$ contains $u, v$, their common neighbors (denote it by $C$), the set of vertices $A \subseteq N(u) \setminus N(v)$ and the set of vertices $B \subseteq N(v) \setminus N(u)$, such that $G[A \cup B]$ is a forest with no isolated points (since each path of length three from $u$ to $v$ is determined by an edge of $G[A \cup B]$). Furthermore, since $G$ is triangle free, all the sets $A, B$ and $C$ are independent, and hence, the faces of $G_0$ are cycles of length 4, 5 or 6.

If a face of $G_0$ is not empty in $G$, then there are no $C_5$ in $G$ that uses vertices both in the interior and the exterior of a face of $G_0$. Indeed such a $C_5$ either uses exactly one vertex inside the face that is adjacent to $u$ and $v$, or exactly two vertices inside the face that form an edge $xy$ with $xu, yv \in E(G)$, since any other possibility gives a triangle in $G$. In either case there would be a $C_5$ containing $u$ and $v$ that is not taken into $G_0$, a contradiction. Therefore, each non-empty (in $G$) face of $G_0$ is a separating cycle in $G$ for which there are no crossing $C_5$'s, unless $V(G_0)=5$ and there is no separating $C_5$ in $G$ (otherwise, we may choose $G_0$ to be a separating $C_5$). However, this implies that every $C_5$ in $G$ is a face and then $\cN(C_5,G)$ is at most the number of faces, which is way smaller than the stated value, and we are done. We consider two cases:

\textbf{Case 1.} $V(G)\setminus V(G_0)=\emptyset$.

Then, $G_0=G$, which means all $5$-cycles of $G$ contain $u$ and $v$. We show that $\cN(C_5,G)=\cN(C_5,G_0) \leq \lf (n-3)/2\rf \cdot \lceil (n-3)/2\rceil$. Paths of length two connecting $u$ and $v$ are corresponding to their common neighbors, that is, there are $|C|$ such paths. Paths of length three connecting $u$ and $v$ must use an edge between $A$ and $B$.  As $G[A \cup B]$ is a forest, there are at most $|A|+|B|-1$ edges between $A$ and $B$ (Note that $G[A \cup B]$ cannot contain an odd cycle as it is bipartite, and cannot contain an even cycle as, together with $u, v$ and a common neighbor of them, it would give a subdivision of $K_{3,3}$ in $G$, contradicting the planarity of $G$). Therefore,
\[\cN(C_5,G)= |C| (|A|+|B|-1) \leq |C|(n-3-|C|) \leq \lf (n-3)/2\rf \cdot \lceil (n-3)/2\rceil,\]
and equality holds only when $|C|$ and $|A|+|B|-1$ are as equal as possible, with $|E(A,B)|=|A|+|B|-1$, which means $G$ is in $\cJ_n$. 

\vskip3mm
\textbf{Case 2.} $V(G)\setminus V(G_0) \neq \emptyset$.

Then, there is a non-empty face of $G_0$, and hence, a separating cycle of $G$ for which there is no crossing $C_5$. As the separating cycle $C_t$ is a face of $G_0$, we have $t$ is $4, 5$ or $6$. Let $L$ be the set of vertices of $G$ in the interior region of $C_t$, with $|L|=l$. Put $G_1:=G[L \cup V(C_t)]$ and $G_2:=G[V(G)\setminus L]$. By the induction hypothesis, we obtain
\begin{align}\label{ineq2}
    \begin{split}
        \cN(C_5,G)&= \cN(C_5, G_1)+\cN(C_5,G_2)\\
        &\leq \left(\frac{(l+t)-3}{2}\right)^2 +  \left(\frac{(n-l)-3}{2}\right)^2\\
        &=\left(\frac{n-3}{2}\right)^2+\frac{t(t-6)+2l^2-2l(n-t)+9}{4}\\
        \end{split}
   \end{align}

The proof is complete if we show that $f(n,t,l):=\frac{t(t-6)+2l^2-2l(n-t)+9}{4} <-1/4$. Observe that the cycle $C_t$ has $l$ vertices in its interior, and it is a separating cycle, then its exterior is also non-empty, so $n\geq l+t+1$. Let us consider the possible values of $t$.

\begin{itemize}
    \item If $t=4$, then using $n\geq l+t+1$, we have $f(n,t,l)\leq \frac{t(t-6)+9-2l}{4}=\frac{1-2l}{4}$. Thus, for $l\geq 2$, we have $f(n,t,l)<-1/4$. If $l=1$, then there is only one vertex inside the separating 4-cycle. We can easily see that in this case $G_1$ contains no $C_5$, and hence $\cN(C_5,G)=\cN(C_5,G_2) \leq \left(\frac{(n-l)-3}{2}\right)^2=\left(\frac{n-4}{2}\right)^2< \lf (n-3)/2\rf \cdot \lceil (n-3)/2\rceil$.

    \item If $t=5$, then the separating 5-cycle is counted twice, which means we must have an extra $-1$. Thus, $f(n,t,l)\leq \frac{t(t-6)+9-2l}{4}-1=-l/2\leq -1/2$.

    \item  If $t=6$, then as the separating 6-cycle is a face in each of the $G_1$ and $G_2$, none of them is isomorphic to $F_{n_i}$ (where $n_i:=|V(G_i)|$, for $i=1,2$). Hence, by the induction hypothesis, each of them contains strictly fewer $C_5$'s than the stated bound in the theorem. So, we get 
    \[\cN(C_5, G_1) \leq \left(\frac{(l+t)-3}{2}\right)^2 -1 \ \ \text{and} \quad \ \ \cN(C_5,G_2)\leq \left(\frac{(n-l)-3}{2}\right)^2-1\]
    Thus, from (\ref{ineq2}), we obtain that $f(n,t,l)=\frac{t(t-6)+2l^2-2l(n-t)+9}{4}-2$. If $l=1$ and there is exactly one vertex in the exterior of $C_t$, then it is easy to see that each of $G_1$ and $G_2$ can have at most two $C_5$'s, then $\cN(C_5,G)\leq 2+2=4 < 6 =\lf (n-3)/2\rf \cdot \lceil (n-3)/2\rceil$, since $n=6+1+1=8$.
    Otherwise, if $l\geq 2$, then using $n\geq l+t+1$, we get $f(n,t,l)=(1-2l)/4\leq -3/4$. Also, if $n\geq l+t+2$, then we have $f(n,t,l)=(1-4l)/4\leq -3/4$.
\end{itemize}
   \end{proof}

For the proof of Theorem \ref{C3C5}  (and Theorem \ref{C3C6} in the next section), we use the concept of \textit{triangular blocks} introduced in \cite{GhGyMPX}. Let us first recall the definition and relevant results from them that we will rely on in the proofs.

\begin{defn} \cite{GhGyMPX} Let $G$ be a plane graph. An edge $e \in E(G)$ is a triangular-block if it is not in any face of length 3 (this is called a \textbf{trivial triangular block}), otherwise we inductively build up the block; start with the subgraph $H:=e$, keep adding to $H$ all faces of length 3 (and their edges) that contain an edge of $H$ until no such is left.
\end{defn}

Let $\cT$ be the set of all triangular blocks of $G$.
It is easy to see that every edge of $G$ is in exactly one triangular block. Then,
\begin{equation}\label{equation2}
    e(G)=\sum_{B \in \cT} e(B). \end{equation}
 
Obviously, every triangle is in at most one triangular block. However, there are possibly triangles in $G$ that are not contained in any triangular block (as a subgraph), let $T$ be the set of all such triangles. Let $\cT_1$ be the set of triangular blocks of $G$ that contain no edges of the triangles in $T$, and $\cT_2:=\cT \setminus \cT_1$. Then 
\begin{equation}\label{equation3}
    \cN(C_3,G)=\sum_{B \in \cT_1} \cN(C_3,B) + \sum_{B \in \cT_2} \cN(C_3,B)+|T|
\end{equation}

The first part of the following proposition was proven in \cite{GyWZh} and the second part was done in \cite{GhGyMPX}.

\begin{prop}\label{tr-blocksinC5C6-free} \cite{GyWZh, GhGyMPX} Let $G$ be a plane graph and $B$ be a triangular block of $G$. Then,
\begin{enumerate}
    \item If $G$ is $C_5$-free, then $B$ has at most four vertices, and it is one of the following graphs:
    \begin{center}
    \begin{tikzpicture}
    \node[vertex] (a) {};
\node [vertex] (b) [right=1.5cm of a] {};
\node [vertex] (c) [above right=1cm and 2cm of b] {};
\node [vertex] (d) [below left=2cm and 0.7cm of c] {};
\node [vertex] (e) [below right=2cm and 0.7cm of c] {};
\node [vertex] (f) [right=4cm of c] {};
\node[vertex] (g) [below left=1cm and 0.7cm of f] {};
\node[vertex] (h) [below right=1cm and 0.7cm of f] {};
\node[vertex] (i) [below=2cm of f] {};
\node [vertex] (j) [right=4cm of f] {};
\node [vertex] (k) [below left=2cm and 1cm of j] {};
\node [vertex] (l) [below right=2cm and 1cm of j] {};
\node [vertex] (m) [below=1.2cm of j] {};

\node[] at (1,-1.5) {$K_2$};
\node[] at (3.5,-1.5) {$K_3$};
\node[] at (7.9,-1.5) {$\Theta_4$};
\node[] at (12,-1.5) {$K_4$};

\path (a) edge (b);
\path (c) edge (d);
\path (c) edge (e);
\path (d) edge (e);
\path (f) edge (g);
\path (f) edge (h);
\path (g) edge (h);
\path (g) edge (i);
\path (i) edge (h);
\path (j) edge (k);
\path (j) edge (l);
\path (j) edge (m);
\path (k) edge (l);
\path (k) edge (m);
\path (l) edge (m);
\end{tikzpicture}
 \end{center}   

    \item If $G$ is $C_6$-free, then $B$ has at most five vertices, and it is either one of the four graphs mentioned above or one of the following:
\begin{center}
    \begin{tikzpicture}
\node [vertex] (c) {};
\node [vertex] (d) [below left=2cm and 0.7cm of c] {};
\node [vertex] (e) [below right=2cm and 0.7cm of c] {};
\node [vertex] (f) [below=0.6cm of c] {};
\node [vertex] (g) [below=1.2cm of c] {};
\node[vertex] (h) [right=3cm of c] {};
\node[vertex] (i) [right=2cm of h] {};
\node[vertex] (j) [below=2cm of h] {};
\node[vertex] (k) [below=2cm of i] {};
\node[vertex] (l) [below right=1cm and 1cm of h] {};
\node[vertex] (m) [right=2cm of i] {};
\node[vertex] (n) [right=2cm of m] {};
\node[vertex] (o) [below=2cm of m] {};
\node[vertex] (p) [below=2cm of n] {};
\node[vertex] (q) [below right=1cm and 0.8cm of n] {};

\node[vertex] (r) [right=3cm of n] {};
\node[vertex] (s) [below left=1cm and 0.6cm of r] {};
\node[vertex] (t) [below=2cm of r] {};
\node[vertex] (u) [below right=1cm and 0.5cm of r] {};
\node[vertex] (v) [right=0.6cm of u] {};

\node[] at (0,-2.7) {$K_5^-$};
\node[] at (4,-2.7) {$B_{5,a}$};
\node[] at (8.5,-2.7) {$B_{5,b}$};
\node[] at (12.7,-2.7) {$B_{5,c}$};

\path (c) edge (d);
\path (c) edge (e);
\path (c) edge (f);
\path (f) edge (g);
\path (d) edge (f);
\path (d) edge (g);
\path (d) edge (e);
\path (e) edge (f);
\path (e) edge (g);
\path (h) edge (i);
\path (h) edge (j);
\path (h) edge (l);
\path (k) edge (i);
\path (k) edge (j);
\path (k) edge (l);
\path (i) edge (l);
\path (j) edge (l);
\path (m) edge (n);
\path (m) edge (o);
\path (o) edge (p);
\path (o) edge (n);
\path (n) edge (p);
\path (n) edge (q);
\path (p) edge (q);
\path (r) edge (s);
\path (r) edge (t);
\path (r) edge (u);
\path (r) edge (v);
\path (t) edge (s);
\path (t) edge (u);
\path (t) edge (v);
\path (u) edge (v);
\end{tikzpicture}
 \end{center}
\end{enumerate}    
\end{prop}

\begin{figure}[ht]
    \centering
    \includegraphics[width=0.9\textwidth]{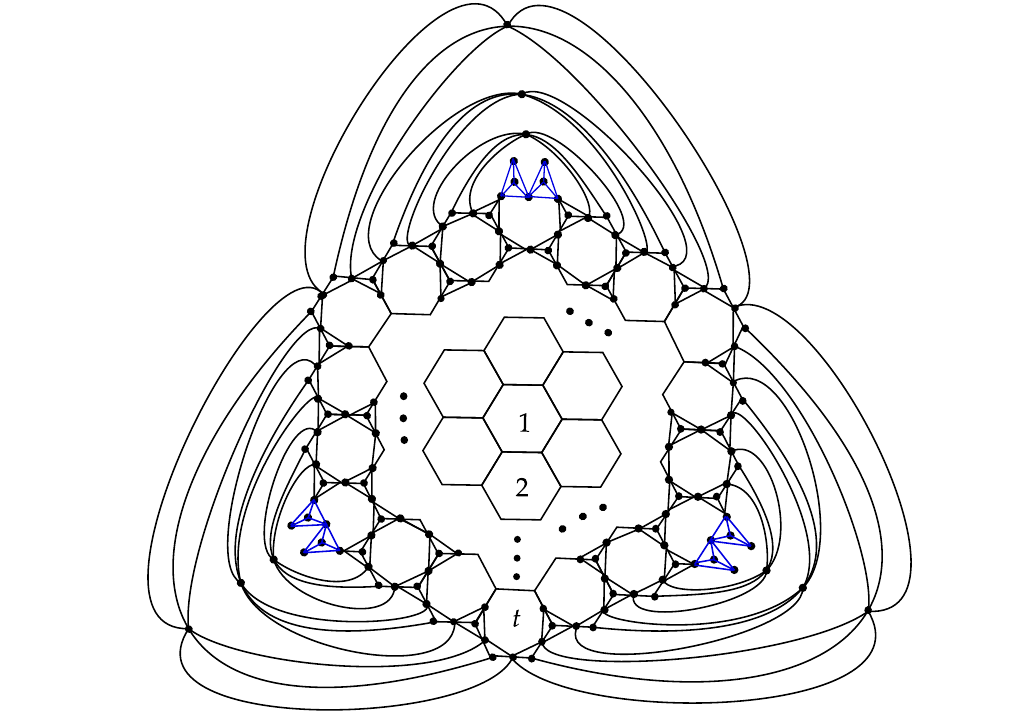}
    \caption{A construction that shows the sharpness in Theorem \ref{C3C5}.}
    \label{C3inC5free}
\end{figure}

\begin{proof}[Proof of Theorem \ref{C3C5}] 
Let $G$ be a $C_5$-free planar graph. Let $\cT, \cT_1, \cT_2$ and $T$ denote the same sets as in (\ref{equation3}). It is easy to see that $\cN(C_3,B) \leq \frac{2}{3}e(B)$ for every triangular block $B$ of $G$, with equality if and only if $B$ is a $K_4$. This immediately, gives 
\begin{equation}\label{ineq4}
    \displaystyle \sum_{B \in \cT_1} \cN(C_3,B) \leq \sum_{B \in \cT_1} \frac{2}{3}e(B). 
\end{equation}

 Observe that no block in $\cT_2$ can be a $K_4$, since it would form a $C_5$ in $G$. Also, for each triangle $t \in T$, either each of its edges are in different blocks or two of its edges are in the same block and the third edge is in a different block. In the former case, at least two out of the three blocks containing the edges of $t$ must be trivial blocks, and $t$ is the only triangle containing them. Thus, for each such $t\in T$, there are at least two trivial blocks in $\cT_2$. In the latter case, the block containing the two edges of $t$ must be a $\Theta_4$ and the other block containing the third edge of $t$ must be a trivial block, i.e. these two blocks form a $K_4$ in $G$ (which is not a triangular block in $G$). Then, we obtain two such triangles in $T$, and two blocks $K_2$ and $\Theta_4$ in $\cT_2$, and no other triangles in $T$ can contain the edges of these blocks. Therefore, computing the number of triangles and edges in $T$ and $\cT_2$, we obtain that 
\begin{equation}\label{ineq5}
    \sum_{B \in \cT_2} \cN(C_3,B)+|T| \leq \sum_{B \in \cT_2} \frac{2}{3}e(B).
\end{equation}

Substituting (\ref{ineq4}) and (\ref{ineq5}) in  (\ref{equation3}), and applying (\ref{equation2}), we obtain

$$\cN(C_3,G)\leq \sum_{B \in \cT}\frac{2}{3}e(B)=\frac{2}{3} \sum_{B \in \cT}e(B)=\frac{2}{3}e(G).$$

Also, it is proved in \cite{GyWZh} (see also \cite{Dowden}) that for every $n\geq 11$, $\ex_\cP(n,C_5)\leq \frac{12n-33}{5}$ , which gives: $$\ex_\cP(n,C_3,C_5) \leq \frac{2}{3} \lf \frac{12n-33}{5}\rf \leq \lf \frac{8n-22}{5}\rf$$

To prove the sharpness of the result, we take the graph Figure \ref{C3inC5free}: we slightly modify the extremal construction for the number of edges (i.e. for $\ex_\cP(n,C_5)$) given in \cite{GyWZh}. In their construction, $n=15t^2-6$, $e=36t^2-21$, and all triangular blocks are $K_4$'s except three of them, which are $\Theta_4$'s. We replace each $\Theta_4$ by two copies of $K_4$ sharing a vertex. Thus, we add 9 vertices and 21 edges, so we get $n=15t^2+3$ and $e=36t^2=\lf\frac{12n-33}{5}\rf$ (This means the new construction is still extremal for the number of edges, too). In the new graph, every triangular block is a $K_4$, and hence, it attains the stated bound on the number of triangles.
\end{proof}

\section{Triangles and 6-cycles}\label{sec4}

Here, we prove Theorems \ref{C3C6} to \ref{P53vtxC3free}. While the proof of Theorem \ref{C3C6} has a simple short proof, which is analogous to the proof of Theorem \ref{C5C3}, the proof of Theorem \ref{C6C3} is more sophisticated.

\begin{proof}[Proof of Theorem \ref{C3C6}]
Let $G$ be a $C_6$-free plane graph, let $\cT, \cT_1, \cT_2$ and $T$ denote the same sets as in (\ref{equation3}). Every triangular-block $B$ of $G$ contains at most $7e(B)/9$ triangles, with equality if and only if $B$ is a $K_5^-$. This implies 
\begin{equation}\label{ineq6}
    \displaystyle \sum_{B \in \cT_1} \cN(C_3,B) \leq \sum_{B \in \cT_1} \frac{7}{9}e(B).
\end{equation}

 Notice that no triangular block in $\cT_2$ can be a $K_5^-$, since this would form a $C_6$ in $G$. 

Assume that each edge of $t$ is in a different triangular block. For $i\in [3]$, let $e_i$ be the edges of $t$, and $B_i\in \cT_2$ be the blocks that contain them such that for each $i$, $e_i \in B_i$. Not that none of the blocks can be a $B_{5,a}$. If for some $i \in[3]$, $B_i \in \{K_4, B_{5,b}, B_{5,c}\}$ (see these blocks in Proposition \ref{tr-blocksinC5C6-free}), then the other two blocks must be trivial, and $t$ is the only triangle that contains them. Thus, for each such triangle in $T$, there are two trivial blocks in $\cT_2$, so the ratio of the number of triangles in these blocks and $t$ to the number of edges in them is less than $7/9$. 

If two of the blocks, say $B_1$ and $B_2$ are not the trivial block, then $B_3$ must be trivial and $B_1$ and $B_2$ are either a $K_3$ or a $\Theta_4$. In this case, $t$ is the only triangle that contains $e_3$, and if any of $e_1$ or $e_2$ is in another triangle $t' \in T$, then the other two edges of $t'$ must be in trivial blocks. Again, the ratio of the number of triangles to the number of edges is less than $7/9$. If one of $B_1$ or $B_2$, say $B_2$ is also trivial, but $e_2$ is in another triangle $t'\in T$, then the other two edges of $t'$ must be in trivial blocks. If further $B_1$ is also trivial, and $e_1$ is contained in another triangle $t'' \in T$, then again the two other edges of $t''$ must be in trivial blocks. Thus, computing such triangles and edges in these blocks, we get to the same conclusion.

Now, assume that two edges of $t$, say $e_1$ and $e_2$ are in the same triangular block $B$ and the third edge $e_3$ is in another block $B'$. Then, $B \in \{\Theta_4, B_{5,a}, B_{5,b}, B_{5,c}\}$. If $B$ is a $\Theta_4$, then $B'$ is a $K_2$, $K_3$ or a $\Theta_4$, and in each case if there is another triangle $t'\in T$ that contains $e_3$, then the other two edges of $t'$ must be in trivial blocks.
Then computing the number of triangles and the number of edges in each possibility we get a ratio of at most $7/9$. Similarly, analyzing the other possibilities of $B\in \{B_{5,a}, B_{5,b}, B_{5,c}\}$, we get the same conclusion.

Therefore, we obtain 
\begin{equation}\label{ineq7}
    \sum_{B \in \cT_2} \cN(C_3,B)+|T| \leq \sum_{B \in \cT_2}\frac{7}{9}e(B).
\end{equation}

Substituting (\ref{ineq6}) and (\ref{ineq7}) in  (\ref{equation3}), and applying (\ref{equation2}), we obtain
\[\cN(C_3,G)\leq \sum_{B \in \cT}\frac{7}{9}e(B)=\frac{7}{9}e(G)\]

It was proved in \cite{GhGyMPX} that for every $n\geq 18$, $\ex_\cP(n,C_6)\leq 5n/2-7$ , which gives: $$\ex_\cP(n,C_3,C_6) \leq \frac{7}{9}\left( \frac{5}{2}n-7 \right) \leq  \frac{35n-98}{18}$$
Sharpness of the bound follows from the construction given in \cite{GhGyMPX} for the number of edges, as every triangular block is a $K_5^-$ that contains $7e(B)/9$ triangles (See Theorem 4 and the construction in Section 2 of \cite{GhGyMPX}).
\end{proof}





\begin{proof}[Proof of Theorem \ref{P5inC3-free}]
    First, suppose that $u$ and $v$ have no common neighbors. We may assume that every edge of $G$ is on a path of length four between them, otherwise we can delete it without decreasing the number of such paths. This implies that adding the edge $uv$ does not violate the planarity of $G$ and since $u$ and $v$ do not have a common neighbor, it does not create a triangle, that is, $G':=G+uv$ is a triangle-free planar graph on $n$ vertices. As each path of length four from $u$ to $v$ in $G$ yields a $C_5$ in $G'$, applying Theorem \ref{C5C3} on $G'$, we have at most $(\frac{n-3}{2})^2 < (\frac{n-1}{2})^2-2$ such paths.

    Now, assume that $u$ and $v$ have $k$ common neighbors. Let $W:=\{w_1, w_2, \ldots, w_k\}$ be the set of common neighbors of $u$ and $v$, ordered in the obvious way, see Figure \ref{kregions}. Then, as $G$ is triangle-free, $W$ is an independent set, and hence, $G[\{u,v\} \cup W]$ is an induced $K_{2,k}$ subgraph of $G$. When drawing this subgraph, it divides the plane into $k$ regions $R_i= uw_ivw_{i+1}$, for $1\leq i\leq k$, where addition in the subscript is modulo $k$. Each path of length four from $u$ to $v$ lies in only one of these regions (may use boundary vertices of the region as well), otherwise it forms a triangle. 

    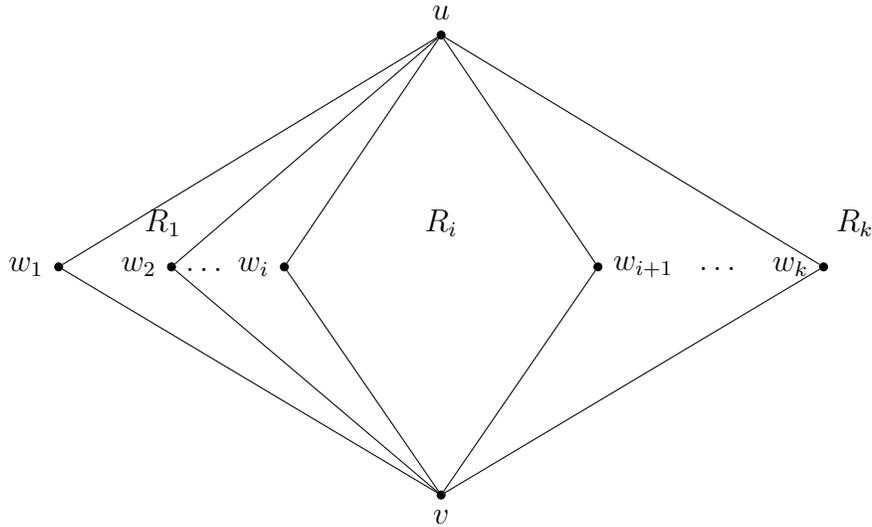
\begin{figure}[ht]
     \begin{center}
    \begin{tikzpicture}
    \node[vertex, label=above: $u$] (a) {};
\node [vertex, label =left: $w_1$] (c) [below left=3cm and 5cm of a] {};
\node [vertex, label =left: $w_2$] (d) [below left=3cm and 3.5cm of a] {};
\node [vertex, label =left: $\ldots \ w_i$] (e) [below left=3cm and 2cm of a] {};
\node [vertex, label =right: $w_{i+1} \ \ \ldots$] (f) [below right=3cm and 2cm of a] {};
\node [vertex, label =left: $w_k$] (g) [below right=3cm and 5cm of a] {};
\node[vertex, label=below: $v$] (b) [below= 6cm of a] {};
\node[] at (0,-2.5) {$R_i$};
\node[] at (-3.7,-2.5) {$R_1$};
\node[] at (5.5,-2.5) {$R_k$};

\path (c) edge (a);
\path (c) edge (b);
\path (d) edge (a);
\path (d) edge (b);
\path (e) edge (a);
\path (e) edge (b);
\path (f) edge (a);
\path (f) edge (b);
\path (g) edge (a);
\path (g) edge (b);
\end{tikzpicture}
 \end{center}
\caption{The common neighbors of $u$ and $v$ form $k$ regions}
        \label{kregions}   
   \end{figure}
    
    Let $M$ be the set of vertices inside the region $R_i$, which is bounded by the 4-cycle $uw_ivw_{i+1}$, with $|M|=m$. Then, together with $u, v, w_i$ and $w_{i+1}$, there are $m+4$ vertices in $R_i$, and $u$ and $v$ have exactly two common neighbors $w_i$ and $w_{i+1}$. Let $X:=\{x_1, x_2, \ldots, x_{m_1}\}$, and $Y:=\{y_1, y_2, \ldots, y_{m_2}\}$ be the neighbours of $u$ and $v$, respectively, in $R_i$ (so $X \cap Y=\{w_i,w_{i+1}\}$). We may assume that every vertex in $R_i$ is on a $P_5(u,v)$, as otherwise deleting it does not reduce the number of such paths. Then, every vertex of $Z:= M \setminus (X \cup Y)$ has neighbors in both $X$ and $Y$ and each vertex of $X$ and $Y$ has neighbors in $Z$. Let $P_5(u,v)_{R_i}$ be a $P_5(u,v)$ that uses only vertices in the region $R_i$ (including $u, v, w_i$ and $w_{i+1}$).

\begin{claim}\label{P5inRi} In the region $R_i$, $\# P_5(u,v)_{R_i} \leq (\frac{m+3}{2})^2-2$.
    \end{claim}
\begin{proof}[Proof of the claim]
 First, we show that the bipartite graph between $Y$ and $Z$ is a forest. Clearly it does not contain an odd cycle, so assume, for a contradiction, that $y_1z_1y_2z_2 \ldots y_rz_r$ is a $C_{2r}$ (an even cycle of length $2r$). If for some $1<j<r$,  $z_1,z_2$ and $z_j$ have a common neighbor $x \in X$, then there is a subdivision of a $K_{3,3}$ with parts $\{v,z_1,z_r\}$ and $\{x,y_1,y_r\}$. If for some $1<j<r$,  $z_1,z_j$ and $z_r$ each has a different neighbor in $X$, then we obtain a subdivision of a $K_{3,3}$ with parts  $\{v,z_1,z_r\}$ and $\{u,y_1,y_r\}$. If for any $1<j<r$, exactly two of the vertices from $\{z_1,z_j,z_r\}$ have a common neighbor in $X$, say $x$, then after contracting the edge $ux$ (to a vertex $u'$), we obtain a subdivision of a $K_{3,3}$ with parts $\{v,z_1,z_r\}$ and $\{u',y_1,y_r\}$. It is well known that contracting an edge in a planar graph results in a planar graph and then should not contain a subdivision of $K_{3,3}$, obtaining a contradiction in any case. Likewise, the bipartite graph between $X$ and $Z$ is a forest. This means $|E(Y,Z)| \leq |Y|+|Z|-1$ and  $|E(X,Z)| \leq |X|+|Z|-1$.

Note that any path $P_5(u,v)$ uses exactly one edge $yz$ between $Y$ and $Z$, and each such edge is in $|N_X(z)|$ paths (except when $y$ is $w_i$ or $w_{i+1}$, for which it is $|N_{X}(z)|-1$), and hence, 
    \begin{equation}\label{eqn-uP5v}
        \#  P_5(u,v)_{R_i} =( \sum_{yz \in E(Y,Z)} |N_X(z)| ) -2.
        \end{equation}
    Let $p=max \{|N_X(z)| \ : \ z \in Z\}$, and let $z_0 \in Z$ have $p$ neighbors in $X$. Then, equation \ref{eqn-uP5v} implies that $\# P_5(u,v)_{R_i} \leq (|Y|+|Z|-1)p-2$. 
   
    We delete all edges between $Y$ and $Z\setminus \{z_0\}$, all edges between $X$ and $Z\setminus \{z_0\}$, and all edges between $Y$ and $X$ (if there is any, these are in no $P_5(u,v)$), and then we add all edges from $Z\setminus \{z_0\}$ to both of $v$ and $z_0$. We also add $w_iz_0$ and $w_{i+1}z_0$, if they are not there already. That is, we move all vertices of $Z$, except $z_0$, to $Y$, obtaining $Y'$, and $Z'=\{z_0\}$. Then, $|Y'|+|Z'|=|Y|+|Z|$, and $|E(Y',Z')|=|Y'|=|Y|+|Z|-1$, with each edge is now in $p$ paths. Hence,  $\#  P_5(u,v)_{R_i} = p|Y'| -2$. Moreover, if $p=|X|$, then any vertex $z \in Z \setminus \{z_0\}$ could have been adjacent to at most one vertex in $X$, which means the edge $yz$ was previously in only one path, and hence the described procedure strictly increases the number of paths. Also, if $p <|X|$, we can now join $z_0$ to all of $X$ and hence strictly increasing the number of paths to $\# P_5(u,v)_{R_i} = |X||Y'|-2\leq |X|(m+3-|X|) -2$ (Note that $|X|+|Y'|=m+3$, since $X \cup Y'=M\setminus \{z_0\}$, and $w_i$ and $w_{i+1}$ are in both of $X$ and $Y'$). This is clearly maximum if $|X|$ and $|Y'|$ are either $\lf(m+3)/2\rf$ or $\lceil (m+3)/2\rceil$, i.e. $|m_1-m_2|\leq 1$, giving that $\#  P_5(u, v)_{R_i} \leq (\frac{m+3}{2})^2-2$.
   \end{proof}
   
Now, we complete the proof of the lemma. If $k=1$, then there is only one region, one common neighbor of $u$ and $v$, and assuming that there are $m$ vertices in this only region, we have $n=m+3$. Then, adapting the proof of the above claim, it is easy to see that $\# \ P_5(u,v) \leq (\frac{m+1}{2})^2-1= (\frac{n-2}{2})^2-1<(\frac{n-1}{2})^2-2$.  

    Now, suppose $k\geq 2$. Assume  two regions $R_i$ and $R_j$ both contain paths of length four from $u$ to $v$, so they are not empty, say they contain $m_i$ and $m_j$ vertices, respectively. Then, applying the claim above, in the two regions together the number of paths is at most 
    \begin{align*}
        \left(\frac{m_i+3}{2}\right)^2+\left(\frac{m_j+3}{2}\right)^2-4&=\frac{m_i^2+m_j^2+6(m_i+m_j)+ 18}{4}-4\\
        &= \frac{(m_i+m_j)^2 -2m_im_j+6(m_i+m_j)+18}{4}-4\\
        &=\left(\frac{m_i+m_j+3}{2}\right)^2-\frac{7+2m_im_j}{4}< \left(\frac{m_i+m_j+3}{2}\right)^2-2.
    \end{align*}

    That is, we obtain more paths if all the $m_i+m_j$ vertices are in one region and the other is empty. Thus, easily by induction, we will have more paths if all the regions except one of them are empty, in which case by Claim \ref{P5inRi} we will have 
    $$ \# \ P_5(u,v) \leq \left(\frac{n-(k+2)+3}{2}\right)^2-2\leq \left(\frac{n-1}{2}\right)^2-2.$$
\end{proof}


Now, we prove three lemmas to prepare the proof of Theorem \ref{C6C3}. In fact, we prove that for sufficiently large $n$ any extremal graph on $n$ vertices has to be isomorphic to $H_n$. Clearly, $H_n$ is a triangle-free planar graph on $n$-vertices.


Let us count the number of 6-cycles in $H_n$. There are $|A||B||C|$ 6-cycles on $V(H_n)\setminus \{z_1,z_2\}$, and each of $z_1$ and $z_2$ is in $|A||B|+|A||C|+|B||C|$ 6-cycles that does not contain the other. Finally, we count those 6-cycles that contain both $z_1$ and $z_2$, any such cycle must contain $u_1,u_2$ and $u_3$, and hence the last vertex comes from one of $A$, $B$ or $C$, but each of them twice depending on which edges we use. Thus, there are $2(|A|+|B|+|C|)$ such 6-cycles. Hence, all together, the number of 6-cycles in $G$ is at most
$$|A||B||C|+2(|A||B|+|A||C|+|B||C|)+2(|A|+|B|+|C|).$$
Note that if $n \equiv 2$ (mod 3), then each of $A, B$ and $C$ contains $(n-5)/3$ vertices, if $n \equiv 1$ (mod 3), then one of them contains $(n-7)/3$ vertices and the others contain $(n-4)/3$ each, and if $n \equiv 0$ (mod 3), then one of them contains $(n-3)/3$ vertices and the other two contain $(n-6)/3$ each. Calculating as above, $H_n$ contains $h(n)$ 6-cycles. Note that among such graphs with various sizes of $A, B$ and $C$, this is the maximum number of 6-cycles.

Also, through Theorem \ref{P5inC3-free}, it is easy to see that the vertices of $H_n$ that are in the largest set among $A, B$ and $C$, are in the fewest number of 6-cycles.  Simply computing that, we see that each vertex of $H_n$ is contained in at least $h_1(n)$ 6-cycles (and there are vertices that are contained in exactly $h_1(n)$ 6-cycles), where $h_1(n)$ is defined as follows.

\begin{align*}\displaystyle
    h_1(n)=\left\{
\begin{array}{ll}
    \left(\frac{n}{3}\right)^2-2 =\frac{n^2}{9}-2, &~ \ \ \text{if} \ n \equiv 0 \ (mod \ 3) \\
     \left(\frac{n-1}{3}\right)\left(\frac{n+2}{3}\right)-2=\frac{n^2}{9}+\frac{n}{9}-\frac{20}{9},  &~ \ \ \text{if} \ n \equiv 1 \ (mod \ 3)\\
      \left(\frac{n+1}{3}\right)^2-2=\frac{n^2}{9}+\frac{2n}{9}-\frac{17}{9},  &~ \ \ \text{if} \ n \equiv 2 \ (mod \ 3)
\end{array} 
\right.
\end{align*}

\vskip3mm
  Note that a vertex of degree one is not contained in any 6-cycle, so we may assume that our graphs have minimum degree at least 2. We fist prove that for $n$ large enough, if an $n$-vertex triangle-free planar graph has the property that each of its vertices is in at least $h_1(n)$ 6-cycles, then it has to be isomorphic to $H_n$. This is done through the following lemmas, following the approach in \cite{savery} and \cite{GhGyJPSZ}.

\begin{lemma}\label{skeleton} Let $n$ be sufficiently large and $0<\gamma <1$ be a constant. Assume that $G$ is a triangle-free planar graph on $n$ vertices such that every vertex of $G$ is contained in at least $n^2/10$ 6-cycles. Then, $G$ contains three vertices $u_1,u_2$ and $u_3$ such that each of $|N(u_1,u_2)|, |N(u_1,u_3)|$ and $|N(u_2,u_3)|$ is at least $\gamma n$.
\end{lemma}
\begin{proof}
Suppose $G$ is a plane graph on $n$ vertices satisfying all the conditions of the lemma. 

\begin{claim}\label{ifK2,an,emptK2,3}
     If $G$ contains two vertices $x$ and $y$ with $|N(x,y)|\geq \alpha n$, where $0<\alpha<1$ is a constant, then $G$ contains a vertex of degree $2$, whose only two neighbors are $x$ and $y$.
\end{claim}
\begin{proof}
    Let $N(x,y)=\{w_1,w_2, \ldots, w_k\}$, then $k \geq \alpha n$. As in the proof of Lemma \ref{P5inC3-free}, they form $k$ regions. For each $1\leq i \leq k$, the region $R_i$ is bounded by the 4-cycle $xw_iyw_{i+1}$, where addition in the subscript is modulo $k$. Since $G$ is triangle-free, $N(x,y)$ is an independent set, and hence, for any $1\leq i \leq k$, the vertex $w_i$ can only be adjacent to $x$, $y$ and vertices inside $R_{i-1}$ or $R_i$. Thus, it suffices to show that there are at least two consecutive empty regions.

    We show that if a region is not empty, then it must contain more than $n^{1/3}$ vertices. On the contrary, assume that a region $R_i$ bounded by $xw_iyw_{i+1}$ is not empty and contains at most $n^{1/3}$ vertices. Let $z$ be a vertex inside $R_i$. We count all the $6$-cycles that contain $z$. If such a $6$-cycle uses vertices only in $R_i$ (including the boundary vertices), then owing to Theorem \ref{asymptC2k}, there are at most $(n^{1/3}+4)^3/9+o((n^{1/3}+4)^2) < n$, as $n$ is large. 

    If a 6-cycle containing $z$ is a crossing one, i.e. uses at least one vertex outside of $R_i$, then it uses at least two vertices on the boundary of $R_i$, and hence it can contain at most three vertices inside $R_i$. 
    
    Assume that the 6-cycle contains only $z$ inside $R_i$. Then, $z$ must be adjacent to at least two of the boundary vertices of $R_i$. As $z$ is not a common neighbor of $x$ and $y$, and since $G$ is triangle-free, the two boundary vertices must be $w_i$ and $w_{i+1}$. Now, without using any of $x$ and $y$ we cannot have a $6$-cycle containing $w_izw_{i+1}$, and if one of $x$ and $y$ are used, we need two vertices outside $R_i$, that form a path of length three from $x$ to $w_i$ or from $x$ to $w_{i+1}$, if $x$ is used (and for $y$ it is similar). A path of length three between $x$ and $w_i$ is determined by an edge between $N(x)$ and $N(w_i)$, since $G$ is triangle-free, $G[N(x) \cup N(w_i)]$ is a forest (similar to the situation in the proof of Theorem \ref{C5C3}), and hence there are at most $n$ choices. Also, if both $x$ and $y$ are used, to complete the 6-cycle we need a common neighbor of $x$ and $y$, for which there are $k-2$ choices. Thus, all together, there are at most $5n$ such 6-cycles.
    
    Assume that the 6-cycle contains $z$ and another vertex, say $z'$, inside $R_i$. There are at most $n^{1/3}$ choices for $z'$. If only two boundary vertices are used, then we need two more vertices outside of $R_i$ that form a path of length three between the two boundary vertices. Like the previous case, there are at most $n$ such paths for each pair of the boundary vertices. Thus, there are at most $6n n^{1/3}=6n^{4/3}$ choices for such a cycle. If the $6$-cycle uses three vertices on the boundary of $R_i$, then together with $z$ and $z'$, it is already 5 vertices, so we need one vertex outside, for which there are at most $n$ choices. Therefore, there are at most $4nn^{1/3}=4n^{4/3}$, as there are 4 possible ways to choose three boundary vertices. Note that it may not be the case that any two or three of the boundary vertices could be chosen, but even this loose upper bound works for us now. In this case, there  are at most $10n^{4/3}$ such 6-cycles.

    Finally, assume that the 6-cycle contains two other vertices $z_1$ and $z_2$ besides $z$ inside the region $R_i$. As it uses at least a vertex outside of $R_i$, it must use exactly two vertices on the boundary of $R_i$ and the outside vertex is a common neighbor of them. Thus, the two boundary vertices cannot be adjacent as this would form a triangle, and hence, they are either $x$ and $y$ or  $w_i$ and $w_{i+1}$. There are no common neighbors of $w_i$ and $w_{i+1}$  outside the region $R_i$ because of the other common neighbors of $x$ and $y$ and planarity of $G$. Thus, the two boundary vertices must be $x$ and $y$ for which there are $k-2$ possible common neighbors outside $R_i$. Then, $x,y,z,z_1$ and $z_2$ form a path of length four from $x$ to $y$ in the region $R_i$ (including boundary vertices), applying Theorem \ref{P5inC3-free}, there are at most $((n^{1/3}+3)/2)^2-2$ such paths, which is less than $n^{5/6}$ as $n$ is large. Even if $z$ is in all those paths, there is at most $kn^{5/6}<n^{11/6}$ 6-cycles in this case.

    Therefore, in total, $z$ is contained in at most $n+5n+10n^{4/3}+n^{11/6} <n^2/10$ 6-cycles, since $n$ is large, obtaining the desired contradiction. Thus, if there are $t$ non-empty regions, they contain more than $tn^{1/3}$ vertices, and hence $tn^{1/3}<n$, which means $t<n^{2/3}$. Since there are $k\geq \alpha n$ regions and $n$ is sufficiently large, there exist consecutive empty regions, which completes the proof. 
\end{proof}

\begin{claim}\label{u1,u2,N(u1,u,2)}
Let $0<\beta <1$ be a constant. If a vertex $v$ has degree at most $5$, then it has two neighbors $w_1$ and $w_2$ for which there is a vertex $w$ such that each of $|N(w,w_1)|$ and $|N(w,w_2)|$ is at least $\beta n$.    
\end{claim}
\begin{proof}
     Every 6-cycle containing $v$ must contain $xvy$ for some pair of neighbors $x$ and $y$ of $v$. There are at most $10$ such pairs and $v$ is in at least $n^2/10$ 6-cycles. Thus, there is a pair $w_1$ and $w_2$ such that $w_1vw_2$ is in at least $n^2/100$ 6-cycles. Each such cycle is determined by a path of length four from $w_1$ to $w_2$, hence there are $n^2/100$ paths $P_5(w_1,w_2)$. Let $P$ be the set of vertices that are on such paths from $w_1$ to $w_2$. Let $X:=N(w_1)\cap P$, $Y:=N(w_2)\cap P$ and $Z:=P\setminus X \cup Y$. Each path uses an edge $yz \in E(Y,Z)$ and a neighbor of $z$ in $X$, and each edge $yz \in E(Y,Z)$ is in at most $|N_X(z)|$ such paths. As in the proof of Theorem \ref{P5inC3-free}, the bipartite graph between $Y$ and $Z$ is acyclic, and hence,  $|E(Y,Z)|\leq |Y|+|Z|-1<n$, and there are at most $\sum_{yz\in E(Y,Z)} |N_X(z)|$ paths of length four from $w_1$ to $w_2$. Let $L:=\{yz\in E(Y,Z): \ |N_X(z)| \geq n/200\}$. We then have $n^2/100 \leq \# P_5(w_1,w_2)\leq |L| n+(n-|L|)n/200$, which implies that $|L| \geq n/199$.

     Since the bipartite graph between $X$ and $Z$ is also acyclic, there are at most $n$ edges between them. Thus, at most 200 vertices in $Z$ can be the endpoints of the edges in $L$, and hence, there is a vertex $w \in Z$, which is the end point of at least $n/39800$ edges in $L$. This means $N(w,w_1)\geq n/199$ and $N(w,w_2) \geq n/39800$. As $n$ is sufficiently large, this completes the proof of the claim.
\end{proof}

Now, we can complete the proof of the lemma. Since $G$ is planar, it contains a vertex of degree at most 5. Then By Claim \ref{u1,u2,N(u1,u,2)}, there are vertices $u_1$ and $u_2$ such that $|N(u_1,u_2)|$ is at least $\gamma n$.
Hence, by Claim \ref{ifK2,an,emptK2,3}, there is a vertex $v_1$ of degree 2, whose only two neighbors are $u_1$ and $u_2$. Applying Claim \ref{u1,u2,N(u1,u,2)} on $v_1$, there is a vertex $u_3$ such that each of $|N(u_1,u_3)|$ and $|N(u_2,u_3)|$ is at least $\gamma n$.
\end{proof}

For a set of vertices $M$, we denote by $P_5(x,y)_M$ a path of length four from $x$ to $y$ that contains some vertices from $M$.

\begin{lemma}\label{P5for3vtx}
    Let $G$ be a triangle-free planar graph, and $u_1v_1u_2v_2u_3v_3$ be an induced 6-cycle in $G$. Suppose $M$ is the set of vertices in the interior of the region bounded by the 6-cycle, with $|M|=m$, such that each vertex in $M$ is adjacent to at most one of $u_1, u_2$ and $u_3$. Then, 
    \[\sum_{i=1}^3 \#P_5(u_i,u_{i+1})_M\leq 3\left(\frac{m+5}{3}\right)^2-3,\]
    where $i+1$ is taken modulo 3.
\end{lemma}

\begin{proof}
First observe that no path $P_5(u_i,u_{i+1})$, for each $i \in [3]$, can use vertices from both the interior and the exterior of the bounding 6-cycle, since it would create a triangle in $G$. Thus, any path $P_5(u_i,u_{i+1})_M$ uses only the vertices of $M$ and those on the 6-cycle.  For each $i \in [3]$, let $X_i$ be the set of neighbors of $u_i$ in $M$ together with it's neighbors $v_i$ and $v_{i+2}$ on the 6-cycle. Since each vertex in $M$ is adjacent to at most one of $u_i$'s, we have $X_i \cap X_{i+1}=\{v_i\}$. Each $P_5(u_i,u_{i+1})$ contains a vertex from $X_i$, a vertex from $X_{i+1}$ and a common neighbor of them, which we call the \textit{middle vertex} of the path.
Let $Z=M\setminus \cup_{i=1}^3X_i$. Then, for each $i\in [3]$, the middle vertex of any $P_5(u_i,u_{i+1})$ is in $Z \cup X_{i+2}$. Accordingly, we distinguish the following cases. 

\vskip3mm
\textbf{Case 1.} For each $i \in [3]$, the middle vertex of any $P_5(u_i,u_{i+1})$ is in $Z$.

For each $i\in [3]$, let $Z_i \subseteq Z$ be the set of middle vertices of all the paths $P_5(u_i,u_{i+1})$. It is easy to see, due to the properties of $G$ being planar and triangle-free, that $|\cap_{i=1}^3Z_i|\leq 1$. Analogous to the proof of Theorem \ref{P5inC3-free}, we have the most number of paths between $u_i$ and $u_{i+1}$ if for each $i \in [3]$, we have $Z_i$ is a singleton. Furthermore, a simple computation gives that we have more paths if $Z_1=Z_2=Z_3=Z=\{z\}$, see Figure \ref{bestforP5for3vtx}. Then, $\sum_{i=1}^3|X_i|\leq m+5$ (since we exclude $z$ from $M$ and each $v_i$ is counted twice), and for each $i \in [3]$, we have $|X_i||X_{i+1}|-1$ paths $P_5(u_i,u_{i+1})$. The term $-1$ is because choosing $v_i$ in both of $X_i$ and $X_{i+1}$ does not give a path of length four from $u_i$ to $u_{i+1}$. This gives 
\[\sum_{i=1}^3 \#P_5(u_i,u_{i+1})_M=\sum_{i=1}^3 \left(|X_i||X_{i+1}| -1\right)\leq 3\left(\frac{m+5}{3}\right)^2-3\]

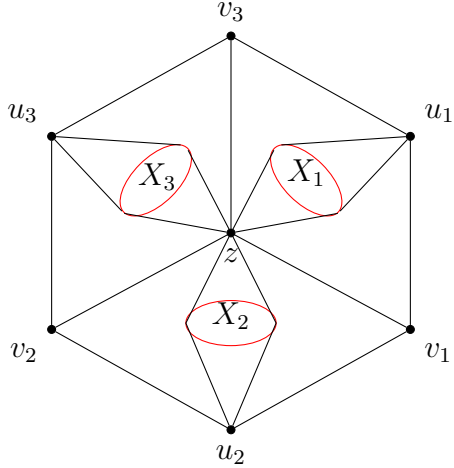
\begin{figure}[ht]
    \begin{center}
        \begin{tikzpicture}
            \node[vertex, label=below:$z$] (z){};
            \node[vertex, label=above right:$u_1$] (u1) [above right=1.2cm and 2.3cm of z] {};
            \node[vertex, label=below right:$v_1$] (v1) [below right=1.2cm and 2.3cm of z] {};
            \node[vertex, label=below:$u_2$] (u2) [below=2.5cm of z] {};
            \node[vertex, label=below left:$v_2$] (v2) [below left=1.2cm and 2.3cm of z] {};
            \node[vertex, label=above left:$u_3$] (u3) [above left=1.2cm and 2.3cm of z] {};
            \node[vertex, label=above:$v_3$] (v3) [above=2.5cm of z] {};

            \draw[red, rotate around={45:(1,0.7)}] (1,0.7) ellipse (0.3cm and 0.6cm);
\draw[red] (0,-1.2) ellipse (0.6cm and 0.3cm);
\draw[rotate around={135:(-1,0.7)},red] (-1,0.7) ellipse (0.3cm and 0.6cm);

\node[] at (1,0.8) {$X_1$};
\node[] at (0,-1.1) {$X_2$};
\node[] at (-0.98,0.75) {$X_3$};

            \draw (u1) -- (v1) -- (u2) -- (v2) -- (u3) -- (v3) -- (u1);
            \draw (1.45,0.29) -- (u1) -- (0.67,1.17);
            \draw (1.42,0.26) -- (z) -- (0.57,1.1);

            \draw (-1.45,0.29) -- (u3) -- (-0.67,1.17);
            \draw (-1.42,0.26) -- (z) -- (-0.57,1.1);

            \draw (z) -- (-0.6,-1.2) -- (u2) -- (0.6,-1.2) -- (z);
            \draw (z) -- (v1);
            \draw (z) -- (v2);
            \draw (z) -- (v3);

        \end{tikzpicture}
    \end{center}
    \caption{The best configuration that achieves the bound in Lemma \ref{P5for3vtx}.}
    \label{bestforP5for3vtx}
\end{figure}

\textbf{Case 2.} For some $i \in [3]$, there is a vertex in $X_{i+2}$, which is the middle vertex of a $P_5(u_{i}, u_{i+1})$.

To ease the notation, without loss of generality, assume there is a vertex $x_3 \in X_3$, which is the middle vertex of a path $P_5(u_1,u_2)$. Since $x_3$ is not adjacent to any  of $u_1$ and $u_2$, it must have neighbors in both of $X_1$ and $X_2$. Let $N_{X_1}(x_3)=Y_1$ and $N_{X_2}(x_3)=Y_2$. Take $x_1 \in Y_1$ such that the region $R_3$ bounded by the cycles $u_1x_1x_3u_3v_3$ does not contain a common neighbor of $u_1$ and $x_3$. Similarly, choose $x_2 \in Y_2$ such that the region $R_2$ bounded by the cycle $u_2x_2x_3u_3v_2$ does not contain a common neighbor of $u_2$ and $x_3$. Let $R_1$ be the region bounded by the cycle $u_1v_1u_2x_2x_3x_1$. For each $i \in [3]$, let $m_i$ be the number of vertices in the interior of $R_i$, so $m_1+m_2+m_3=m-3$. See Figure \ref{middleinX3}.

\begin{figure}[ht]
    \begin{center}
        \begin{tikzpicture}
            \node[] (z){};
            \node[vertex, label=above right:$u_1$] (u1) [above right=1.2cm and 2.3cm of z] {};
            \node[vertex, label=below right:$v_1$] (v1) [below right=1.2cm and 2.3cm of z] {};
            \node[vertex, label=below:$u_2$] (u2) [below=2.5cm of z] {};
            \node[vertex, label=below left:$v_2$] (v2) [below left=1.2cm and 2.3cm of z] {};
            \node[vertex, label=above left:$u_3$] (u3) [above left=1.2cm and 2.3cm of z] {};
            \node[vertex, label=above:$v_3$] (v3) [above=2.5cm of z] {};

            \node[vertex, label=above:$x_3$] (x3) [above left=0.7cm and 1cm of z] {};

            \node[vertex, label=above: {\small $x_1$}] (x1) [above right= 1.2cm and 0.8cm of z]{};
            \node[vertex] (x12) [below right= 0.2cm and 0.1cm of x1]{};
            \node[vertex] (x13) [below right= 0.2cm and 0.1cm of x12]{};

            \node[vertex, label=left:{\small $x_2$}] (x2) [below left=1.3cm and 0.6cm of z] {};
            \node[vertex] (x22) [right=0.2cm of x2] {};
             \node[vertex] (x23) [right=0.2cm of x22] {};

             \node[] at (0,2.2) {$R_3$};
             \node[] at (1.4,-0.8) {$R_1$};
             \node[] at (-2,-0.8) {$R_2$};

            \draw (u1) -- (v1) -- (u2) -- (v2) -- (u3) -- (v3) -- (u1);
            \draw (u3) -- (x3);
            \draw (u1) -- (x1) -- (x3) -- (x2) -- (u2);
            \draw (u1) -- (x12) -- (x3) -- (x22) -- (u2);
            \draw (u1) -- (x13) -- (x3) -- (x23) -- (u2);
        \end{tikzpicture}
    \end{center}
    \caption{$x_3\in X_3$ is the a middle vertex of some $P_5(u_1,u_2)$.}
    \label{middleinX3}
\end{figure}
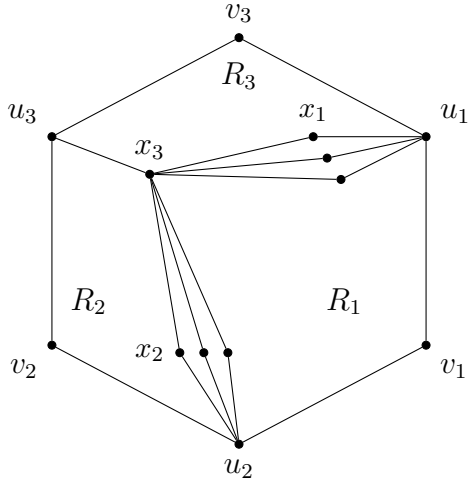

Note that, since $G$ is planar and triangle-free, none of the $x_i$'s can be any of the $v_i$'s, and no other vertices of $X_3$ can be the middle vertex of a $P_5(u_1,u_2)$. Also, for the same reasons, if $Y_1$ contains only $x_1$, then it is possible that $x_1$ is the middle vertex of a path $P_5(u_2,u_3)$, but then no vertex in $X_2$ can be the middle vertex of a path $P_5(u_1,u_3)$. That is, in each $X_i$ at most one vertex can be the middle vertex, and at most two of the $X_i$'s can contain a middle vertex. We now have the following possibilities:

\vskip3mm
\textbf{Subcase 2.1} Each of $Y_1$ and $Y_2$ contains at least two vertices. Thus, none of $x_1$ and $x_2$ can be the middle vertex of the considered paths. Let us count the paths according the regions $R_1, R_2$ and $R_3$.

No vertex inside $R_3$ can be on a path $P_5(u_1,u_2)$ or $P_5(u_2,u_3)$, and by Theorem \ref{P5inC3-free} there are at most $\left(\frac{m_3+4}{2}\right)^2-2$ paths $P_5(u_1,u_3)$ (note that there are $m_3$ vertices inside $R_3$ and 5 vertices on its boundary). Similarly, vertices inside $R_2$ cannot be on paths $P_5(u_1,u_2)$ and $P_5(u_1,u_3)$, and there are at most $\left(\frac{m_2+4}{2}\right)^2-2$ paths $P_5(u_1,u_3)$. 

Let us consider $R_1$. Again there at most $\left(\frac{m_1+5}{2}\right)^2-2$ paths $P_5(u_1,u_2)$. If a path $P_5(u_1,u_3)$ contains vertices inside $R_1$, then it must contain the edge $x_3u_3$ and a path $P_4(u_1,x_3)$, i.e. a path of length three from $u_1$ to $x_3$. Likewise, a $P_5(u_2,u_3)$ that uses vertices inside $R_1$ is determined by a $P_4(u_2,x_3)$. Let $L_1$ and $L_2$ be the sets of vertices inside $R_1$ that are on paths $P_4(u_1,x_3)$ and $P_4(u_2,x_3)$. We claim that $|L_1 \cap L_2|\leq 1$. To see this, assume a vertex $x'\in X_1$ (a neighbor of $u_1$) is on a path $P_4(u_2,x_3)$. Then $x'$ must be adjacent to $x_3$ and must have a neighbor with $u_2$, since $x'$ is not adjacent to $u_2$. Thus, $x'$ is a common neighbor of $u_1$ and $x_3$, and hence, it cannot be on a $P_4(u_1,x_3)$, as otherwise this would form a triangle. Therefore, no vertex in $X_1$ (inside $R_1$) can be in $L_1 \cap L_2$. Similarly, no vertex in $X_2$ inside $R_1$ can be in $L_1 \cap L_2$. Now, suppose a vertex $w$ inside $R_1$ is in $L_1 \cap L_2$. Then $w$ is not adjacent to any of $u_1$ and $u_2$ and it is on a $P_4(u_1,x_3)$ and a $P_4(u_2,x_3)$. Then, it must be adjacent to $x_3$ and have neighbors with both $u_1$ and $u_2$, which then separates the region $R_1$ into three other regions that can easily be seen that no other vertices can be common to both of $L_1$ and $L_2$. As mentioned other times earlier, a $P_4(u_1,x_3)$ is determined by an edge in $L_1$, which is an acyclic bipartite graph between the neighbors of $u_1$ and $x_3$, and hence, there are at most $|L_1|-1$ such paths. Similarly, there are at most $|L_2|-1$ paths $P_4(u_2,x_3)$. Also, on the boundary, $v_1$ can be counted in $L_1 \cap L_2$. Thus, there are at most $m_1+4$ paths of length three from $u_1$ or $u_2$ to $x_3$, which means there are at most $m_1+4$ paths $P_5(u_1,u_3)$ or $P_5(u_2,u_3)$. Thus, we obtain

\begin{align*}
    \sum_{i=1}^3 \#P_5(u_i,u_{i+1})_M&= \left(\frac{m_3+4}{2}\right)^2-2+ \left(\frac{m_2+4}{2}\right)^2-2+\left(\frac{m_1+5}{2}\right)^2-2+m_1+4\\
    &= \left(\frac{m_3+4}{2}\right)^2+ \left(\frac{m_2+4}{2}\right)^2+\left(\frac{m_1+5}{2}\right)^2+m_1-2\\
    &< 3\left(\frac{m+5}{3}\right)^2-3.
    \end{align*}

\textbf{Subcase 2.2} At least one of $Y_1$ or $Y_2$ contains only one vertex. Then, it is possible for $x_1$ or $x_2$ to be the middle vertex of a considered path. Without loss of generality, assume $x_1$ is such, which means $Y_1=\{x_1\}$. Then, $x_1$ has neighbors in $X_2$, which are not in $Y_2$ (as otherwise with $x_3$, they form a triangle) and may have other neighbors besides $x_3$ in $X_3$, which must be in the region $R_3$. 

If $x_1$ has no other neighbors besides $x_3$ in $X_3$, then this situation is already covered in Subcase 2.1. And if $x_3$ has other neighbors besides $x_3$ in $X_3$, then the situation is the same as the previous subcase with the role of $x_1$ and $x_3$ being swapped.
\end{proof}

\begin{lemma}\label{extgraphisHn} Let $G$ be an $n$-vertex triangle-free planar graph. If every vertex of $G$ is contained in at least $h_1(n)$ 6-cycles, then $G$ is isomorphic to $H_n$, if $n$ is sufficiently large.  
\end{lemma}

\begin{proof}
    Let $G$ be an $n$-vertex triangle-free planar graph, where $n$ is sufficiently large. Assume every vertex of $G$ is contained in at least $h_1(n)$ 6-cycles. By Lemma \ref{skeleton}, there are three vertices $u_1, u_2$ and $u_3$ such that for every $i \in [3]$, $|N(u_i,u_{i+1})|\geq \alpha n$, where $0<\alpha <1$ is a constant. The addition in the subscript, here and in the rest of the proof, is taken modulo 3. Let $|N(u_i,u_{i+1})|=k_i$, for each $i \in [3]$. Let $G_0$ be the induced subgraph of $G$ on $\{u_1, u_2, u_3\}\cup N(u_1,u_2) \cup N(u_2,u_3) \cup N(u_1,u_3)$. 

    Consider $G[\{u_1, u_2\} \cup N(u_1,u_2)]$. This divides the plane into $k_1$ regions, each bounded by a 4-cycle $u_1xu_2y$, for some $x,y \in N(u_1,u_2)$. Then, $u_3$ and its common neighbors with $u_1$ and with $u_2$ are in one of these regions, without loss of generality, say in the ``first" region, as shown in Figure \ref{skeletonfig}

    Note that $|N(u_1,u_2,u_3)|\leq 2$, as otherwise $G$ would contain a $K_{3,3}$. If $|N(u_1,u_2,u_3)|=2$, then all the faces of $G_0$ are bounded by 4-cycles, if $|N(u_1,u_2,u_3)| =1$, then there is one face bounded by a 6-cycle, if $|N(u_1,u_2,u_3)|=0$, there are two faces bounded by 6-cycles.  See Figure \ref{skeletonfig}. We first show that $V(G)\setminus V(G_0)=\emptyset$.

    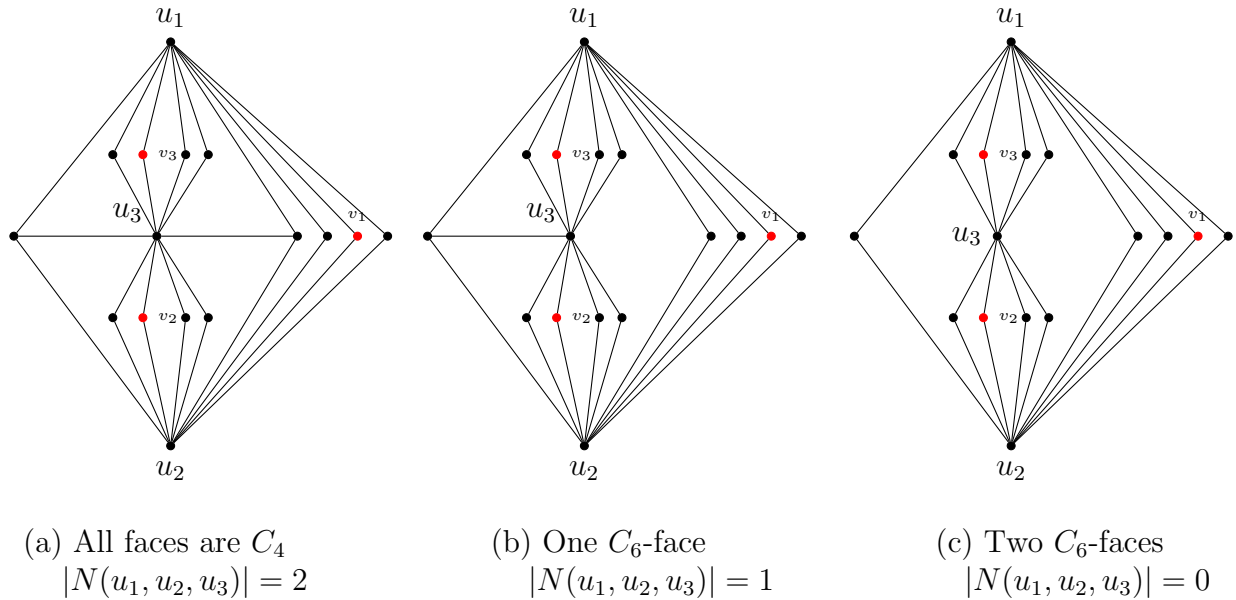
\begin{figure}[ht]
        \begin{center}
            \begin{tikzpicture}
                \node[] (O) {};
                \node[vertex, label=above: $u_1$] (u1) [above left=2.5cm and 7cm of O] {};
                \node[vertex, label=above left: $u_3$] (u3) [below left=2.5cm and 0.1cm of u1] {};
                \node[vertex, label=below: $u_2$] (u2) [below left=2.5cm and 7cm of O] {};
                
                \node[vertex, label=above: $u_1$] (u12) [above left=2.5cm and 1.5cm of O] {};
                \node[vertex, label=above left: $u_3$] (u32) [below left=2.5cm and 0.1cm of u12] {};
                \node[vertex, label=below: $u_2$] (u22) [below left=2.5cm and 1.5cm of O] {};
                
                \node[vertex, label=above: $u_1$] (u13) [above right=2.5cm and 3.8cm of O] {};
                \node[vertex, label=left: $u_3$] (u33) [below left=2.5cm and 0.1cm of u13] {};
                \node[vertex, label=below: $u_2$] (u23) [below right=2.5cm and 3.8cm of O] {};

                \node[vertex] (w1) [below left=2.5cm and 2cm of u1] {};
                \node[vertex] (w2) [below right=2.5cm and 1.6cm of u1] {};\node[vertex] (w3) [below right=2.5cm and 2cm of u1] {};
                \node[vertex, red, label=above: {\tiny $v_1$}] (v1) [below right=2.5cm and 2.4cm of u1] {};
                \node[vertex] (wk) [below right=2.5cm and 2.8cm of u1] {};

                \node[vertex] (x1) [below left=1cm and 0.5cm of u3] {};
                \node[vertex, red, label=right: {\tiny $v_2$}] (v2) [below left=1cm and 0.1cm of u3] {};
                \node[vertex] (x2) [below right=1cm and 0.3cm of u3] {};
                \node[vertex] (xk) [below right=1cm and 0.6cm of u3] {};

                \node[vertex] (y1) [above left=1cm and 0.5cm of u3] {};
                \node[vertex, red, label=right: {\tiny $v_3$}] (v3) [above left=1cm and 0.1cm of u3] {};
                \node[vertex] (y2) [above right=1cm and 0.3cm of u3] {};
                \node[vertex] (yk) [above right=1cm and 0.6cm of u3] {};

                \node[vertex] (w12) [below left=2.5cm and 2cm of u12] {};
                \node[vertex] (w22) [below right=2.5cm and 1.6cm of u12] {};
                \node[vertex] (w32) [below right=2.5cm and 2cm of u12] {};
                \node[vertex, red, label=above: {\tiny $v_1$}] (v12) [below right=2.5cm and 2.4cm of u12] {};
                \node[vertex] (wk2) [below right=2.5cm and 2.8cm of u12] {};

                  \node[vertex] (x12) [below left=1cm and 0.5cm of u32] {};
                \node[vertex, red, label=right: {\tiny $v_2$}] (v22) [below left=1cm and 0.1cm of u32] {};
                \node[vertex] (x22) [below right=1cm and 0.3cm of u32] {};
                \node[vertex] (xk2) [below right=1cm and 0.6cm of u32] {};

                \node[vertex] (y12) [above left=1cm and 0.5cm of u32] {};
                \node[vertex, red, label=right: {\tiny $v_3$}] (v32) [above left=1cm and 0.1cm of u32] {};
                \node[vertex] (y22) [above right=1cm and 0.3cm of u32] {};
                \node[vertex] (yk2) [above right=1cm and 0.6cm of u32] {};

                \node[vertex] (w13) [below left=2.5cm and 2cm of u13] {};
                \node[vertex] (w23) [below right=2.5cm and 1.6cm of u13] {};\node[vertex] (w33) [below right=2.5cm and 2cm of u13] {};
                \node[vertex, red, label=above: {\tiny $v_1$}] (v13) [below right=2.5cm and 2.4cm of u13] {};
                \node[vertex] (wk3) [below right=2.5cm and 2.8cm of u13] {};

                \node[vertex] (x13) [below left=1cm and 0.5cm of u33] {};
                \node[vertex, red, label=right: {\tiny $v_2$}] (v23) [below left=1cm and 0.1cm of u33] {};
                \node[vertex] (x23) [below right=1cm and 0.3cm of u33] {};
                \node[vertex] (xk3) [below right=1cm and 0.6cm of u33] {};

                \node[vertex] (y13) [above left=1cm and 0.5cm of u33] {};
                \node[vertex, red, label=right: {\tiny $v_3$}] (v33) [above left=1cm and 0.1cm of u33] {};
                \node[vertex] (y23) [above right=1cm and 0.3cm of u33] {};
                \node[vertex] (yk3) [above right=1cm and 0.6cm of u33] {};

                \draw (u1) -- (w1) -- (u2);
                \draw (u1) -- (w2) -- (u2);
                \draw (u1) -- (w3) -- (u2);
                \draw (u1) -- (v1) -- (u2);
                \draw (u1) -- (wk) -- (u2);

               \draw (u2) -- (x1) -- (u3) -- (y1) -- (u1);
               \draw (u2) -- (v2) -- (u3) -- (v3) -- (u1);
               \draw (u2) -- (x2) -- (u3) -- (y2) -- (u1);
               \draw (u2) -- (xk) -- (u3) -- (yk) -- (u1);
               
                \draw (w1) -- (u3) -- (w2);

                 \draw (u12) -- (w12) -- (u22);
                \draw (u12) -- (w22) -- (u22);
                \draw (u12) -- (w32) -- (u22);
                \draw (u12) -- (v12) -- (u22);
                \draw (u12) -- (wk2) -- (u22);

                 \draw (u22) -- (x12) -- (u32) -- (y12) -- (u12);
               \draw (u22) -- (v22) -- (u32) -- (v32) -- (u12);
               \draw (u22) -- (x22) -- (u32) -- (y22) -- (u12);
               \draw (u22) -- (xk2) -- (u32) -- (yk2) -- (u12);

                \draw (w12) -- (u32);

                \draw (u13) -- (w13) -- (u23);
                \draw (u13) -- (w23) -- (u23);
                \draw (u13) -- (w33) -- (u23);
                \draw (u13) -- (v13) -- (u23);
                \draw (u13) -- (wk3) -- (u23);

               \draw (u23) -- (x13) -- (u33) -- (y13) -- (u13);
               \draw (u23) -- (v23) -- (u33) -- (v33) -- (u13);
               \draw (u23) -- (x23) -- (u33) -- (y23) -- (u13);
               \draw (u23) -- (xk3) -- (u33) -- (yk3) -- (u13);

               \node[] at (-7.4,-4) {(a) All faces are $C_4$};
               \node[] at (-1.5,-4) {(b) One $C_6$-face};
               \node[] at (4.5,-4) {(c) Two $C_6$-faces};
               
               \node[] at (-7,-4.5) {$\abs{N(u_1,u_2,u_3)}=2$};
               \node[] at (-0.8,-4.5) {$\abs{N(u_1,u_2,u_3)}=1$};
               \node[] at (5,-4.5) {$\abs{N(u_1,u_2,u_3)}=0$};
            \end{tikzpicture}
        \end{center}
        \caption{The possible cases of the subgrapb $G_0$.}
        \label{skeletonfig}
    \end{figure}
    
    From Claim \ref{ifK2,an,emptK2,3} in the proof of Lemma \ref{skeleton}, we obtain that for each $i\in [3]$, there is a vertex $v_i\in N(u_i,u_{i+1})$ such that the only neighbors of $v_i$ are $u_{i}$ and $u_{i+1}$. For each $i\in [3]$, we count all the 6-cycles that contain $v_i$. Then, by averaging, there is a vertex $v\in \{v_1, v_2, v_3\}$ that is contained in at most $\frac{1}{3} \sum_{i=1}^3 \# C_6(v_i)$ (recall that $C_6(v_i)$ denotes a 6-cycle that contains the vertex $v_i$).
    If a 6-cycle contains $v_i$, then it must contain the path $u_iv_iu_{i+1}$ and a $P_5(u_i,u_{i+1})$. Also, any $P_5(u_i,u_{i+1})$ gives a 6-cycle containing $v_i$. Hence, $\# C_6(v_i)=\# P_5(u_i,u_{i+1})$, which means
    \[\sum_{i=1}^3 \# C_6(v_i)= \sum_{i=1}^3 \# P_5(u_i,u_{i+1})\]
    
    Assume a face of $G_0$ that is bounded by a 4-cycle $u_ixu_{i+1}y$, where $x,y \in N(u_i,u_{i+1})$, for some $i\in [3]$, is not empty in $G$ (i.e. there are vertices of $G$ in its interior). Let $T$ be the set of vertices in its interior with $|T|=t$. Then as in the proof of Claim \ref{ifK2,an,emptK2,3}, we must have $t\geq n^{1/3}$. Observe that if none of $x$ and $y$ is adjacent to $u_{i+2}$, then no vertex of $T$ can be on a $P_5(u_i,u_{i+2})$ or a $P_5(u_{i+1}, u_{i+2})$. Also, at most one of $x$ and $y$ can be adjacent to $u_{i+2}$ (and this is possible only for two such faces), say $x$ is adjacent to $u_{i+2}$. Then, any $P_5(u_i, u_{i+2})$ that uses vertices from $T$ must use a path of length three from $u_i$ to $x$ and the edge $xu_{i+2}$. Then there are at most $t$ such paths. Similarly there are at most $t$ paths $P_5(u_{i+1}, u_{i+2})$ that uses vertices from $T$.  
    By Theorem \ref{P5inC3-free}, $T$ contributes to at most $(t+3)^2/4-2=t^2/4+3t/2+1/4$ paths $P_5(u_i,u_{i+1})$. Thus, by moving $t/2$ of the vertices to each of $N(u_i, u_{i+2})$ and $N(u_{i+1}, u_{i+2})$, we lose at most $2t$ paths of length four from $u_i$ to $u_{i+2}$ and from $u_{i+1}$ to $u_{i+2}$, and at most $t^2/4+3t/2+1/4$ paths of length four between $u_i$ and $u_{i+1}$, while we gain $t^2/4+(k_{i+1}+ k_{i+2})t/2$ paths between $u_{i}$ and $u_{i+1}$, $k_it/2$ paths between $u_i$ and $u_{i+2}$ and $k_it/2$ paths between $u_{i+1}$ and $u_{i+2}$. 
    
    Thus, the total increase of $\sum_{i=1}^3 \# P_5(u_i,u_{i+1})$ is at least 
    \[t^2/4+(k_{i+1}+k_{i+2})t/2+k_it - \left(t^2/4+3t/2+1/4 +2t\right)\geq (k_{i+1}+k_{i+2})t/2 \geq \alpha n^{4/3}.\] 
    
    Therefore, we may assume that all the regions that are bounded by 4-cycles of $G_0$ are empty. We now show that the possible 6-faces of $G_0$ (i.e are bounded by 6-cycles of $G_0$) are empty in $G$, too. Recall that we can have at most two such faces, so suppose $M$ and $L$ are the sets of vertices in their interiors, with $|M|=m$ and $|L|=l$. 
    
    For each $i\in [3]$, a $P_5(u_i,u_{i+1})$ can not contain vertices from both of $M$ and $L$, and hence it either contains no vertex from $M \cup L$ or it contains some vertices from $M$ or it contains some vertices from $L$. Then, owing to Lemma \ref{P5for3vtx}, and assuming $k_1+k_2+k_3=k$, we obtain
    \begin{align*}
        \sum_{i=1}^3 \# C_6(v_i)&= \sum_{i=1}^3 \left( k_ik_{i+1}+ \#P_5(u_i,u_{i+1})_M+ \#P_5(u_i,u_{i+1})_L \right)\\
        & \leq k_1k_2+k_2k_3+k_1k_3+ 3\left(\frac{m+5}{3}\right)^2-3 + 3\left(\frac{l+5}{3}\right)^2-3\\
        &\leq 3\left(\frac{k}{3}\right)^2+3\left(\frac{m+5}{3}\right)^2-3+ 3\left(\frac{l+5}{3}\right)^2-3
    \end{align*}
    
    Then, by averaging, there is a vertex $v \in \{v_1, v_2, v_3\}$ such that
    \begin{align*}
    \# C_6(v)&\leq 1/3 \sum_{i=1}^3 \# C_6(v_i)=
    \left(\frac{k}{3}\right)^2+\left(\frac{m+5}{3}\right)^2+\left(\frac{l+5}{3}\right)^2-2\\
    &=\frac{k^2}{9}+ \frac{m^2+10m+25}{9}+\frac{l^2+10l+25}{9}-2\\
    &=\frac{k^2+m^2+l^2}{9}+ \frac{10m+25}{9}+\frac{10l+25}{9}-2\\
    &=\frac{(k+m+l)^2-2(km+kl+ml)}{9}+\frac{10m+25}{9}+\frac{10l+25}{9}-2\\
    &<\frac{(n-1)^2}{9}-2 < h_1(n),
    \end{align*}
where the last inequality is because $k+m+l \leq n-1$. Indeed, if both of $m$ and $l$ are non-zero, then $k+m+l=n-3$, as they count all the vertices of $G$ except the $u_i$'s. Also, if one of them, say $l$ is zero, then we can have a vertex in $N(u_1,u_2,u_3)$, which will be counted three times in $k_i$'s, but still $k+m=n-1$, and we reach the same conclusion. This contradiction proves that both of $M$ and $L$ must be empty, which gives $V(G_0)=V(G)$ as desired.  

Now, to show that $G$ is isomorphic to $H_n$, we must have all $k_i$'s are as equal as possible and $|N(u_1,u_2,u_3)|=2$. Since all faces are empty, then we can add appropriate chords in the two faces that are bounded by 6-cycles, and make $|N(u_1,u_2,u_3)|=2$, increasing the number of paths of length four between each pair of the $u_i$'s. Therefore, we may assume that $|N(u_1,u_2,u_3)|=2$, and hence $k_1+k_2+k_3=n+1$. Then, for each $i\in [3]$, we have the number of $P_5(u_i,u_{i+1})$, and hence the number of $6$-cycles that contain $v_i$, is exactly $k_{i+1}k_{i+2}-2$.

Assume $n\equiv 0$ (mod 3), which means $n+1\equiv 1$ (mod 3). If for some $i\in [3]$, we have $k_i \geq n/3+1$, then $k_{i+1}+k_{i+2} \leq 2n/3$. Then, $\#C_6(v_i)=\#P_5(u_i,u_{i+1})= k_{i+1}k_{i+2}-2\leq (n/3)^2-2=h_1(n)$, and equality holds only when each of $k_{i+1}$ and $k_{i+2}$ is $n/3$, which means we must also have equality in $k_i \geq n/3+1$, which implies that $G$ is isomorphic to $H_n$.

Assume $n\equiv 1$ (mod 3), which means $n+1\equiv 2$ (mod 3). If for some $i\in [3]$, we have $k_i \geq (n+2)/3$, then $k_{i+1}+k_{i+2} \leq (2n+1)/3$. Note that, since $n\equiv 1$ (mod 3), we have $(2n+1)/3$ is an odd integer, we can then have $\#C_6(v_i)=\#P_5(u_i,u_{i+1})= k_{i+1}k_{i+2}-2\leq ((n+2)/3)\cdot ((n-1)/3)-2=h_1(n)$. Again, equality holds only when one of $k_{i+1}$ and $k_{i+2}$ is $(n+2)/3$, and the other is $(n-1)/3$, which means we must also have equality in $k_i \geq (n+2)/3$, implying that that $G$ is isomorphic to $H_n$.

Finally, if $n\equiv 2$ (mod 3), we can similarly get the same conclusion.
 \end{proof}

We are now ready to prove Theorem \ref{C6C3}.

\begin{proof}[Proof of Theorem \ref{C6C3}]
    Let $n$ be sufficiently large, and $G_n$ be an extremal graph on $n$-vertices. By Lemma \ref{extgraphisHn}, we have either $G_n$ contains a vertex in fewer than $h_1(n)$ 6-cycles or it is isomorphic to $H_n$. So in any case $G_n$ contains a vertex $v$ in at most $h_1(n)$ 6-cycles. Then, by deleting $v$ we obtain a graph $G_{n-1}$ on $n-1$ vertices that contains at least $\ex_\cP(n,C_6,C_3)-h_1(n)$ 6-cycles, which implies $\ex_\cP(n-1,C_6,C_3)\geq \ex_\cP(n,C_6,C_3)-h_1(n)$, i.e.
    \[ \ex_\cP(n,C_6,C_3)-\ex_\cP(n-1,C_6,C_3) \leq h_1(n).\]

On the other hand, we have $h(n)-h(n-1)=h_1(n)$. Therefore, if for if for infinitely many values of $n$ an extremal graph $G_n$ contains a vertex in strictly fewer than $h_1(n)$ 6-cycles, we obtain $\ex_\cP(n,C_6,C_3)<h(n)$, which is a contradiction, since $H_n$ is a triangle-free planar graph containing $h(n)$ 6-cycles. Thus, there is an $n_0$ such that if $G$ is an extremal graph on $n>n_0$ vertices, then every vertex of $G$ is contained in at least $h_1(n)$ 6-cycles. Then, by Lemma \ref{extgraphisHn}, they are isomorphic to $H_n$.
\end{proof}

Finally we present the proof of Theorem \ref{P53vtxC3free}. We will be rather sketchy, as many of the arguments have already been used in the previous proofs. Again, the addition that appears in the subscripts is taken modulo 3. 

\begin{proof}[Proof of Theorem \ref{P53vtxC3free}]
    Let $u_1,u_2$ and $u_3$ be three distinct vertices of $G$. Let $|N(u_i,u_{i+1})|=k_i$, for each $1\leq i \leq 3$, and let $k=k_1+k_2+k_3$. First, assume that each $k_i$ is non-zero, in  particular this means  $\{u_1,u_2,u_3\}$ is independent, since $G$ is triangle-free. Consider $G_0:=G[\cup_{i=1}^3N(u_i,u_{i+1}) \cup  \{u_1,u_2,u_3\}]$. As $G$ is planar, $|N(u_1,u_2,u_3)| \leq 2$. If $|N(u_1,u_2,u_3)|=2$, then we have all the faces of $G_0$ are bounded by 4-cycles, and similar to the proof of Lemma \ref{extgraphisHn}, we obtain more paths in total if all the faces are empty and $G=G_0$. Consequently, we have 
    \[\sum_{i=1}^3 \#P_5(u_i,u_{i+1})=\sum_{i=1}^3 (k_ik_{i+1}-2)\leq 3\left(\frac{n+1}{3}\right)^2-6\]
    where the inequality is because $k_1+k_2+k_3=n+1$ (since they count all the vertices except the $u_i$'s and the two common neighbors are each counted 2 times). The equality holds if each $k_i$ is exactly $\frac{n+1}{3}$, which happens when $n\equiv 2$ (mod 3). (Note that in this case $G \cong H_n$.)

    If $|N(u_1,u_2,u_3)|=1$, then $G_0$ has a face bounded by a 6-cycle. If $k_i\geq 2$, for all $i\in [3]$, then this 6-cycle is an induced one. Assume $M$ is the set of vertices in its interior with $|M|=m$. Applying Lemma \ref{P5for3vtx}, we have 
    \begin{align*}
        \sum_{i=1}^3 \#P_5(u_i,u_{i+1})&=\sum_{i=1}^3 (k_ik_{i+1}-1)+\sum_{i=1}^3 \#P_5(u_i,u_{i+1})_M\\
        &\leq \sum_{i=1}^3 (k_ik_{i+1}-1)+3\left(\frac{m+5}{3}\right)^2-3 \leq 3\left(\frac{k}{3}\right)^2-3+3\left(\frac{m+5}{3}\right)^2-3
        \end{align*}
    Clearly, this is maximum if $k$ or $m$ is the smallest possible. Note that $k\geq 6$ (and counts the vertex in $N(u_1,u_2,u_3)$ three times), and hence, $m$ can be as large as $n-7$. A simple computation shows that we have the maximum possible value if $m=0$ (and $k=n-1$), which is then still much smaller than the stated bound in the theorem.

    Assume for some $i\in [3]$, we have $k_i=1$. Without loss of generality assume $k_2=1$, then $N(u_2,u_3)=N(u_1,u_2,u_3):=\{v\}$. The $C_6$-face of $G_0$ is not an induced one then, and the edge $u_1v$ can be in paths of length four from $u_1$ to $u_2$ or to $u_3$. Any such path must use a path of length three from $v$ to $u_2$ or to $u_3$. Similar to the Subcase 2.2 in the proof of Lemma \ref{P5for3vtx} we can show that there are at most $m+1$ such paths. Thus, there are at most $m+1$ paths of length four from $u_1$ to $u_2$ or $u_3$ that contain the edge $u_1v$. Therefore we have 
\begin{align*}
        \sum_{i=1}^3 \#P_5(u_i,u_{i+1})&\leq \sum_{i=1}^3 (k_ik_{i+1}-1)+\sum_{i=1}^3 \#P_5(u_i,u_{i+1})_M\\
        &\leq \sum_{i=1}^3 (k_ik_{i+1}-1)+3\left(\frac{m+5}{3}\right)^2-3+m+1\\
        &=k_1+k_3+k_2k_3-3+3\left(\frac{m+5}{3}\right)^2-3+m+1 \quad [\text{since} \ k_2=1]\\
        &\leq (k-1)+\left(\frac{k-1}{2}\right)^2-3+3\left(\frac{m+5}{3}\right)^2-3+m+1\\
        &\leq \frac{n^2+n+1}{3}-8,
        \end{align*}
where the last inequality holds because we have the maximum value in case $k=3$ (the smallest possible) and $m=n-6$.
    
    If $|N(u_1,u_2,u_3)|=0$, Then $G_0$ has two such $C_6$-faces, each of which is an induced 6-cycle. Again, applying Lemma \ref{P5for3vtx} to the interior of the $C_6$-faces and counting all the paths as before, we can deduce that there are strictly fewer paths than the stated bound in the theorem.
    
Finally, Suppose for some $i \in [3]$, we have $k_i=0$. Then, again for each $i\in [3]$, a $P_5(u_i,u_{i+1})$ uses a vertex in $N(u_i)$, a vertex in $N(u_{i+1})$ and a middle vertex. One can repeat the argument in the proof of Lemma \ref{P5for3vtx} to see that the number of all paths is maximum if there is only one middle vertex and the rest of the vertices are evenly distributed into the neighbors of $u_1, u_2$ and $u_3$. Then going through the computations, one sees that there are strictly fewer paths $P_5$ joining all the pairs of $u_i$'s than the stated bound in the theorem.
    
\end{proof}
\section{Concluding Remarks}
 The maximum number of an even cycle in triangle-free planar graphs is asymptotically the same as in planar graphs in general. This is clearly true for $C_4$ and $C_6$ as we have shown. For any even cycle $C_{2k}$, the following  construction attains the best asymptotic value, which is given in Theorem \ref{asymptC2k}. 
 
 Take an even cycle $C_{2k}$, and blow up every other vertex in a balanced way such that all blow-up subsets intersect in two vertices, call such graphs $\cG_{\mathrm{even}}$ (Note that $H_n$ is $\cG_{\mathrm{even}}$ for $k=3$). We conjecture that $\cG_{\mathrm{even}}$ is the unique extremal graph for $\ex_\cP(n,C_{2k},C_3)$. 

 However, in case of the odd cycles, we have already seen that $\ex_\cP(n,C_5,C_3)$ is not asymptotically the same as $\ex_\cP(n,C_5,\emptyset)$. In this case, we conjecture the following (similar to $\cJ_n$) construction $\cG_{\mathrm{odd}}$ to be the extremal graphs. Take a cycle $C_{2k+1}=v_1v_2v_3 \ldots v_{2k}v_{2k+1}$ and join two vertices $z_1$ and $z_2$ (one from inside and the other form the outside of $C_{2k+1}$) to each vertex $v_{2i+1}$ for $i=1,2, \ldots, k-1$. Then, blow-up each vertex $v_{2i}$, for $i=1,2, \ldots, k-1$, to independent sets of vertices, and replace the edge $v_{2k}v_{2k+1}$ by a tree with color classes $A$ and $B$, such that the size of each of the blown-up sets and $|A \cup B|-1$ are as equal as possible.

\begin{conj}\label{conj}
    Let $n$ be sufficiently large and $k\geq 2$. Then,
    \begin{enumerate}
        \item $\ex_\cP(n,C_{2k}, C_3)=\cN(C_{2k},\cG_{\mathrm{even}})$, and $\cG_{\mathrm{even}}$ is the unique extremal graph.

        \item $\ex_\cP(n,C_{2k+1}, C_3)=\cN(C_{2k},G_n)$, for a $G_n \in \cG_{\mathrm{odd}}$, and $\cG_{\mathrm{odd}}$ is the set of all extremal graphs.
    \end{enumerate}
\end{conj}

Note that the Theorems \ref{C4C3} and \ref{C6C3} show that the first part of this conjecture holds for $k=2,3$ and Theorem \ref{C5C3} yields part 2 for $k=2$.

\section*{Acknowledgments} Research of Gy\H{o}ri was supported by the National Research, Development and Innovation Office - NKFIH under the grants K 132696  and SNN 135643. Research of Hama Karim was supported by the National Research, Development and Innovation Office - NKFIH under the grant FK 132060.

\end{document}